\documentclass{conm-p-l}

\usepackage{amsfonts}
\usepackage{amssymb,latexsym,amsmath, amscd}

\copyrightinfo{2000}{American Mathematical Society}

\usepackage{amssymb,latexsym,amsmath}
\begin{document}

\def\no{\noindent}

\newtheorem{dfn}{Definition}[section]
\newtheorem{add}[dfn]{Addendum}
\newtheorem{rem}[dfn]{Remark}
\newtheorem{thm}[dfn]{Theorem}
\newtheorem{mthm}[dfn]{Main Theorem}
\newtheorem{lem}[dfn]{Lemma}
\newtheorem{sublem}[dfn]{Sublemma}
\newtheorem{prop}[dfn]{Proposition}
\newtheorem{prob}[dfn]{Problem}
\newtheorem{ass}[dfn]{Assumption}
\newtheorem{classprob}[dfn]{Classical Problem}
\newtheorem{eigenvaluesofasum}[dfn]{Eigenvalues of a sum Problem}
\newtheorem{cor}[dfn]{Corollary}
\newtheorem{conj}[dfn]{Conjecture}
\newtheorem{ex}[dfn]{Example}
\newtheorem{ques}[dfn]{Question}
\newtheorem{techques}[dfn]{Technical Question}
\newtheorem{crit}[dfn]{Criterion}
\newtheorem{listof}[dfn]{List of Properties}
\newtheorem{conv}[dfn]{Convention}
\newtheorem{cons}[dfn]{Consequence}
\newtheorem{defn}[dfn]{Definition}
\newtheorem{fact}[dfn]{Fact}
\newtheorem{obs}[dfn]{Observation}
\newtheorem{warn}[dfn]{Warning}
\newtheorem{stabcrit}[dfn]{Stability Criterion}
\newenvironment{theorema}{\noindent{\bf Theorem A.}\em}{}
\newenvironment{theoremareal}{\noindent{\bf Theorem A (Real case).}\em}{}
\newenvironment{theoremb}{\noindent{\bf Theorem B.}\em}{}
\newenvironment{theoremcsymplectic}{\noindent{\bf Theorem C
        (Symplectic case).}\em}{}
\newenvironment{theoremd}{\noindent{\bf Theorem D.}\em}{}

\let\lhd\vartriangleleft
\def\proof{\par\medskip\noindent{\it Proof. }}
\def\sketch{\par\medskip\noindent{\it Sketch of proof. }}
\def\lra{\longrightarrow}
\def\Lra{\Longrightarrow}
\def\ra{\rightarrow}
\def\Ra{\Rightarrow}
\def\CR{\curvearrowright}
\def\lh{\hookleftarrow}
\def\half{\frac{1}{2}}
\def\C{{\Bbb C}}
\def\R{{\Bbb R}}
\def\E{{\Bbb E}}
\def\H{{\Bbb H}}
\def\Z{{\Bbb Z}}
\def\K{{\mathcal K}}
\def\d{{\mathcal D}}
\def\P{{\Bbb P}}
\def\Q{{\Bbb Q}}
\def\S{{\Bbb S}}
\def\B{{\Bbb B}}
\def\N{{\Bbb N}}
\def\F{{\mathcal F}}
\def\eps{\epsilon}
\def\al{\alpha}
\def\be{\beta}
\def\ga{\gamma}
\def\Ga{\Gamma}
\def\de{\delta}
\def\De{\Delta}
\def\Si{\Sigma}
\def\si{\sigma}
\def\L{{\mathcal L}}
\def\la{\lambda}
\def\La{\Lambda}
\def\Om{\Omega}
\def\om{\omega}
\def\D{\partial}
\def\hook{\hookrightarrow}
\def\embed{\hookrightarrow}
\def\8{\infty}
\def\<{\langle}
\def\>{\rangle}
\def\e{\sim}
\def\BE{\begin{equation}}
\def\EE{\end{equation}}
\def\LL{{\mathcal L}}
\def\geo{\partial_{\infty}}
\def\tits{\partial_{Tits}}
\def\tangle{\angle_{Tits}}
\def\ov{\overrightarrow}
\def\ol{\overline}
\def\grad{\mathop{\hbox{grad}}}
\def\8{\infty}
\newcommand{\restr}{\mbox{\Large \(|\)\normalsize}}
\def\oo{{\cal O}}
\def\too{\tilde{\cal O}}
\def\3{\ss}

\def \pol{{\mathcal P}_n(\mathfrak g)}

%%% New Macros
\newcommand{\bbCPm}{\bbC\Bbb{P}^m}

\def\check{\centerline{\bf CHECK THIS!}}
\def\gap{\centerline{\bf GAP!}}

% temporary macro
\def\goth{\mathfrak}

\title[The generalized triangle inequalities]{The generalized triangle
inequalities for  rank $3$ symmetric spaces of noncompact type.}
\author{Shrawan Kumar}

\address{Department of Mathematics, University of North Carolina,
Chapel Hill, NC 27599-3250, USA.}

\email{shrawan@email.unc.edu}

\author{Bernhard Leeb}
\address{Mathematisches Institut, Universit\"at M\"unchen, Theresienstrasse
39 D-80333 M\"unchen, Germany }
%\email{}

\author{John Millson}

\address{Department of Mathematics, University of Maryland, College
Park, MD 20742, USA.}

\email{jjm@math.umd.edu}

\thanks{S. Kumar was partially supported by NSF grant number
DMS 0070679 and J. Millson was partially supported by NSF grant
number DMS 0104006}

\begin{abstract}
We compute the generalized triangle inequalities explicitly for
all rank $3$ symmetric spaces of noncompact type. For SL$(4,\C)$
there are $50$ inequalities none of them redundant by
\cite{KnutsonTaoWoodward}. For both Sp$(6,\C)$ and Spin$(7,\C)$
there are $135$ inequalities of which $24$ are trivially redundant
in the sense that they follow from the inequalities defining the
Weyl chamber $\Delta$. There are $9$ more redundant inequalities
for each of these two groups. One interesting feature is that
these inequalities do not occur for the other system (and
consequently must be redundant because the two polyhedral cones
are the same by Theorem \ref{transfer}). The two equal  polyhedral
cones $D_3(B_3)=D_3(C_3)$ have precisely $102$ facets and $51$
generators (edges).
\end{abstract}

\dedicatory{Dedicated to Robert Greene on his sixtieth birthday.}

\maketitle

\section{Introduction}

Let $G$ be a connected  semisimple real Lie group with no compact
factors and finite center and Lie algebra $\mathfrak{g}$, $K$ be a
maximal compact subgroup and $X= G/K$ be the associated symmetric
space. By a symmetric space of noncompact type we will mean a
symmetric space $G/K$ where $G$ is as above. We let $\mathfrak{k}$
denote the Lie algebra of $K$ and let $\mathfrak{g} = \mathfrak{k}
\oplus \mathfrak{p}$ be the Cartan decomposition. Let
$\mathfrak{a}$ be a Cartan subspace in $\mathfrak{p}$ (i.e., a
maximal subalgebra in $\mathfrak{p}$ which is
 necessarily
abelian). Let $A$ be the real points of the split torus in $G$
corresponding to $\mathfrak{a}$. Choose an ordering of  the
restricted roots and let $\Delta \subset \mathfrak{a}$ be the
corresponding (closed) Weyl chamber. Let $A_{\Delta} $ be the
image of $\Delta$ under the exponential map $\exp:\mathfrak{g} \to
G$. Let $o$ be the point in $X$ that is fixed by $K$. We will
refer to $o$ as the basepoint for $X$.  We will need the following
theorem, the Cartan decomposition for the group $G$, see
\cite{Helgason}, Theorem 1.1, pg. 402.

\begin{thm}\label{Cartandecomposition}
We have
$$G = KA_{\Delta}K.$$
Moreover, for any $g \in G$, the intersection of the double coset
$KgK$ with $A_{\Delta}$ consists of a single point to be denoted
$a(g)$.
\end{thm}

Suppose $\overline{x_1x_2}$ is the oriented geodesic segment in
$X$ joining the point $x_1$ to the point $x_2$. Then there exists
a unique element $g \in G$ which sends $x_1$ to $o$ and $x_2$ to
$y = \exp(\delta)$ where $\delta \in \Delta$. Note that the point
$\delta$ is uniquely determined by $\overline{x_1x_2}$. We define
a map $\sigma$ from $G$-orbits of oriented geodesic segments to
$\Delta$ by

$$\sigma(\overline{x_1x_2}) = \delta.$$

Clearly we have the following consequence of the Cartan
decomposition.

\begin{lem}
The map $\sigma$ gives rise to a one-to-one correspondence between
$G$--orbits of oriented geodesic segments in $X$ and the points of
$\Delta$.
\end{lem}

\begin{rem} In the real-rank $1$ case $\sigma(\overline{x_1x_2})$ is just
the length of the geodesic segment $\overline{x_1x_2}$. In general
we will call $\sigma(\overline{x_1x_2})$ the $\Delta$--length of
$\overline{x_1x_2}$ or the $\Delta$--distance between $x_1$ and
$x_2$. We will write $d_{\Delta}(x_1,x_2) =
\sigma(\overline{x_1x_2})$.
\end{rem}

 We note the formula
 $$d_{\Delta}(x_1,x_2) = \log a(g_1^{-1}g_2) \ \text{where} \ x_1 = g_1K,
x_2 =g_2K.$$

\begin{rem}
The delta distance is {\em symmetric} in the sense that
$$d_{\Delta}(x_1,x_2) = - w_0 d_{\Delta}(x_2,x_1),$$
where $w_0$ is the longest element in the restricted Weyl group.
However, the naive triangle inequality
$$d_{\Delta}(x_1,x_3) \leq d_{\Delta}(x_1,x_2) + d_{\Delta}(x_2,x_3)$$
does not hold \cite{KapovichLeebMillson2}. Here the order is the
one defined by the cone $\Delta$. The naive triangle inequality
has to be replaced by the inequalities below.
\end{rem}

We have the fundamental problem of finding the {\em generalized
triangle inequalities}  for $X$, precisely we have

\begin{prob}
Give conditions on a triple $(v_1,v_2,v_3) \in \Delta^3$ that are
necessary and sufficient in order that there exist a triangle in
$X$ with vertices $x_1,x_2,x_3$ such that $d_{\Delta}(x_1,x_2) =
v_1$ , $d_{\Delta}(x_2,x_3) = v_2$ and $d_{\Delta}(x_3,x_1) =
v_3$.
\end{prob}

We now describe a system of linear inequalities on $\Delta ^3$
which will give the required necessary and sufficient conditions.
These inequalities will be called the {\em generalized triangle
inequalities}.

\begin{rem}
We will include the inequalities defining $\Delta^3$ in
$\mathfrak{a}^3$ in the generalized triangle inequalities.
\end{rem}

We will need the following:
\begin{dfn}
Suppose that $W$ and $W^{\prime}$ are Coxeter groups acting on $V$
and $V^{\prime}$ respectively by their natural reflection
representations. We define a monomorphism of Coxeter systems to be
a pair $(f,\phi)$  where  $f:V \to V^{\prime}$ is an isometric
embedding and $\phi$ is a monomorphism $W \to W'$ satisfying
$f(wx) =\phi(w)f(x).$
\end{dfn}

First we reduce to the case in which $G$ is complex by the
following:

\begin{thm}[\cite{LeebMillson,KapovichLeebMillson1}] \label{transfer}
\label{cone}

1. The set $D_3(X)\subset \De^3$ of triples $(v_1,v_2,v_3)$ for
which a triangle in
 the Problem 1.5 (for $X$) exists,
is a polyhedral cone.

2. $D_3(X)$ depends only on the spherical Coxeter complex
associated to $X$. More precisely, a monomorphism  $(f,\phi)$  of
the Coxeter system $(\mathfrak{a},W)$ to the Coxeter system
$(\mathfrak{a}^{\prime},W^{\prime})$ induces an affine embedding
$D_3(X)\to D_3(X')$. In particular, if $f$ and $\phi$ are also
surjective, then the map $D_3(X)\to D_3(X')$ is an affine
isomorphism.
\end{thm}

{\em Thus given $G$ as above we can replace $G$ by any complex
semisimple group of the same rank as $G$ whose Weyl group
coincides with the restricted Weyl group $W$ of $G$.} Thus it
suffices to find the generalized triangle inequalities for the
case in which $G$ is complex.

{\it So from now on we will assume that $G$ is a connected complex
semisimple group.} We will accordingly often rewrite  $D_3(X)$ as
$D_3(R)$ where $R$ is the reduced root system associated to the
restricted root system of $X$.

The system of inequalities for $D_3(R)$ breaks up into
rank$(\mathfrak{g})$ subsystems $*_P$ where $P$ is a standard
maximal parabolic subgroup. The subsystem $*_P$ is controlled by
the Schubert calculus in the generalized Grassmannian $G/P$ in the
sense that there is one inequality $*_{w_1,w_2,w_3}$ for each
triple of elements $w_1,w_2,w_3 \in W^P$ such that
$$X^P_{w_1}\cdot X^P_{w_2} \cdot X^P_{w_3} = [pt]$$
in  $H_*(G/P)$. Here $W^P$ is the set of minimal length coset
representatives for the set of cosets $W/W_P$ ($W_P$ being the
Weyl group of $P$), $\cdot$ is the intersection product and
$X^P_w, w \in W^P$, is the Schubert class in $G/P$. To describe
the inequality $*_{w_1,w_2,w_3}$,  let $P$ be the standard maximal
parabolic corresponding to a fundamental weight $\lambda$. Then
the action of $W$ on the weight lattice  of $\mathfrak{g}$ induces
a one-to-one correspondence $f:W^P \to W\lambda$. Thus  we may
reparametrize the Schubert classes in $G/P$ by elements of
$\mathfrak{a}^*$. We let $\lambda_i = f(w_i), i = 1,2,3$. We will
sometimes denote the Schubert cycle $X^P_w$ as $X_{f(w)}$. Then
the inequality $*_{w_1,w_2,w_3}$ is given by
$$\lambda_1(v_1) + \lambda_2(v_2) + \lambda_3(v_3) \leq 0,
\ (v_1,v_2,v_3) \in \Delta^3.$$

\begin{rem} The key point is that there is a basis of the algebra
$H_*(G/P)$ parametrized by certain linear functionals on
$\mathfrak{a}$.
\end{rem}

In order to give an accurate account of the history of work on the
problem, we first need to describe the corresponding problem for
the {\em infinitesimal} symmetric space $\mathfrak{p}$. The $Ad K$
orbits in $\mathfrak{p}$ are again parametrized by the points of
the cone $\Delta$. Thus given a triple of side-lengths
$(v_1,v_2,v_3) \in \Delta^3$ as above we can look for a triple $e
= (e_1,e_2,e_3) \in \mathfrak{p}^3$ such that $e_1 \in Ad K \cdot
v_1$, \ $e_2 \in Ad K \cdot v_2$, \ $e_3 \in Ad K \cdot v_3$ and
$$e_1 + e_2 + e_3 = 0.$$
We let $D_3(\mathfrak{p})$ be the subset of $(v_1,v_2,v_3) \in
\Delta^3$ such that there is a solution to the above problem. The
following theorem (for all connected semisimple groups with no
compact factors) was proved in \cite{KapovichLeebMillson1}
(however see below).

\begin{thm} \label{infinitesimaltriangleinequalities}
$D_3(\mathfrak{p}) = D_3(X).$
\end{thm}

Many people contributed to finding the generalized triangle
inequalities. The reader is urged to consult \cite{Fulton} and
\cite{KapovichLeebMillson2} for a more complete account. In fact
the inequalities were computed for the {\em infinitesimal}
symmetric space $\mathfrak{p}$ for $G$ complex by
\cite{BerensteinSjamaar} and \cite{LeebMillson} and for general
$G$ (real or complex) in \cite{LeebMillson}. An intense interest
in recent years in the generalized triangle inequalities was
triggered by the  paper of Klyachko, \cite{Klyachko1}. Klyachko
proved the generalized triangle inequalities for
$GL(m,\mathbb{C})$ in the infinitesimal symmetric space case and a
refinement was obtained by Belkale \cite{Belkale}. In a second
paper, Klyachko, \cite{Klyachko2} , proved the equality
$D_3(\mathfrak{p}) = D_3(X)$ for certain symmetric spaces of
complex simple groups (including $SL(m)$). In \cite{AMW},
Alekseev, Meinrenken and Woodward proved the corresponding
equality  for all complex simple $G$.

The above system of inequalities is not so explicit; in
particular,  the polyhedral cone $D_3(R)$ is not well understood.
Thus it is important to compute explicit examples. We will see
that the  theorem of Knutson,Tao and Woodward
\cite{KnutsonTaoWoodward} that the inequalities are irredundant
for $GL(m)$ is highly exceptional.

\begin{dfn}
The inequality  $*_{w_1,w_2,w_3}$ will be said to be {\em
trivially redundant} if it is a consequence of the inequalities
defining the cone $\Delta$ in $\mathfrak{a}$.
\end{dfn}

In the case $G$ has rank one, the generalized triangle
inequalities are precisely the ordinary triangle inequalities. The
rank two examples were worked out in \cite{LeebMillson}. For the
cases of $B_2 = C_2$ and $G_2$ the above system of inequalities
was not minimal. For the case of $B_2$ there was one trivially
redundant inequality. Once it was removed the remaining
inequalities were irredundant. For the case of $G_2$ there were no
trivially redundant inequalities but three redundant inequalities.
The point of this paper is to work out all the rank three
examples. We will see that for $C_3$ and $B_3$ there are $24$
trivially redundant inequalities. It is rather surprising that
there are only $9$ more redundant inequalities. Furthermore as
explained in the abstract these redundancies were easy to find
since the redundant inequalities occurred for one system and not
for the other. Thus we might say that only ``obvious''
redundancies occurred. Another interesting consequence is that
although the polyhedral cones $D_3(B_3)$ and $D_3(C_3)$ are
isomorphic (i.e.,  there is an affine isomorphism from one to the
other)
 by Theorem \ref{transfer}, the systems of inequalities are
different.

In Chapter 3, see Theorem \ref{realization}, we have given a
self-contained account of how one uses the Demazure--B.G.G.
operators to realize the duals of the Schubert homology classes in
the Borel model. Our calculations in this paper of the products of
Schubert classes aare based on this realization.

This paper is dedicated to Robert Greene on the occasion of his
sixtieth birthday. The third author takes great pleasure  in
acknowledging many helpful conversations and some hard-fought
tennis matches over the years. This paper depends on
\cite{LeebMillson} and the two papers \cite{KapovichLeebMillson1}
and \cite{KapovichLeebMillson2}. We have used the computer program
Porta written by Thomas Christof and Andreas L\"obel. Finally we
would like to thank George Stantchev for finding the computer
program Porta and for much help and advice in implementing it.

\section{Schubert cycles in generalized flag varieties}

We continue to assume that $G$ is a connected complex semisimple
algebraic group. For convenience, we further assume that $G$ is
simply-connected. We fix a Borel subgroup $B$ of $G$. Let
$\mathfrak{b}$, resp. $\mathfrak{g}$, be the Lie algebra of $B$,
resp. $G$. We also fix a Cartan subalgebra $\mathfrak{h} \subset
\mathfrak{b}$. The choice of  $(\mathfrak{h},
 \mathfrak{b})$ determines the set $\Pi \subset \mathfrak{h}^*$ of positive
roots and thus the set $\Phi =\{\alpha_1, \cdots,
\alpha_l\}\subset \Pi$ of simple roots and also the fundamental
weights $\{\omega_1, \cdots, \omega_l\}, \,l$ being the rank of
$G$.  We consider the real form $\mathfrak{a}$ of
 $\mathfrak{h}$ which is the real span of the simple coroots
$\{\alpha_1^\vee, \cdots, \alpha_l^\vee\}$. Then $\Pi \subset
\mathfrak{a}^*$ and also $\omega_i\in  \mathfrak{a}^*$. Any
algebraic subroup $P$ of $G$ containing $B$ is called a {\it
standard parabolic subgroup}. Let $\Delta \subset \mathfrak{a}^*$
be the cone on the fundamental weights. By definition, the {\it
dominant weights} are the elements of $\oplus_{i=1}^l \,\Bbb
Z_+\omega_i$.

 To each dominant weight
$\lambda $ we define the {\it associated standard parabolic
subgroup} $P_{\lambda}$ to be the (connected) subgroup of $G$ with
Lie algebra spanned by $\mathfrak{b}$ together with the root
vectors $X_{- \alpha}$ (corresponding to the root $-\alpha$) such
that $\lambda(H_{\alpha})=0$. Here $\alpha$ runs through the
positive roots $\Pi$ of $\mathfrak{g}$ and $ H_{ \alpha} \in
\mathfrak{a}$ is the corresponding coroot which is equal to
$\frac{2\alpha}{\langle\alpha, \alpha\rangle}$ under the Killing
form.
 We let $S$ be the set of (simple) reflections in the root
hyperplanes defined by the simple roots and let $W\subset
\text{Aut}\, \mathfrak h$ be the Weyl group generated by $S$. Then
$W$ can canonically be identified with the group $N(T)/T$, where
$T$ is the maximal torus of $G$ with Lie algebra $\mathfrak h$ and
$N(T)$ is its normalizer in $G$.

For any standard parabolic subgroup $P$ of $G$,  let $W_P \subset
W$ be the subgroup consisting of those $w$ that have
representatives in $P$ and $S_P = S \cap W_P$. We let $\ell$ be
the length function on $W$. Each coset $wW_P$ in $W/W_P$  has a
unique representative of minimum length \cite{Hiller}, Ch. I,
Corollary 5.4. We will denote this element by $w^*=w^*_P$. The set
of such representatives will be denoted by $W^P$.  We have the
following criterion for the minimum length element in the coset
$wW_P$, see \cite{Hiller}, Ch. I, Corollary 5.4.

\begin{lem}\label{minimumlengthrepresentative}
$w^* \in wW_P$ is the minimum length representative if and only if
$$\ell(w^* s)=l (w^*) + 1, \ \text{for all} \ s \in S_P.$$
\end{lem}

We will also need the following result, see \cite{Hiller}, Ch. I,
Theorem 5.3.

\begin{lem} \label{lengthsadd}
Suppose that $w \in W^P$ and $v \in W_P$. Then
$$\ell(wv) = \ell(w) + \ell(v).$$
\end{lem}

We now recall the Bruhat decomposition for $G$, see
\cite{Helgason}, Theorem 1.4, pg. 403.

\begin{thm}\label{Bruhatdecomposition} For any standard parabolic subgroup
$P$ of $G$,
\begin{align}
G = & \sqcup_{w \in W} BwB\,, \\
G = & \sqcup_{w \in W^P} BwP.
\end{align}
\end{thm}

As a consequence of $(2)$  the generalized flag variety $G/P$ is
the disjoint union of the subsets $\{C_w^P:= BwP/P\}_{w\in W^P}$.
The subset $C_w^P$ is biregular isomorphic  to the affine space
$\Bbb C^{\ell(w)}$ and is called a {\em Schubert cell}. The
closure of $C_w^P$, to be denoted $X_w^P$, is an algebraic
subvariety of $G/P$ and is called a {\em Schubert variety}. We
will use $[X_w^P]$ to denote the integral homology class in
$H_*(G/P)$ carried by $X_w^P$ but we will often abuse notation and
use the same symbol $X_w^P$ for the variety and its class in
homology. We will use $PD(X_w^P)$ to denote the cohomology class
(of complementary degree to the degree of $X_w^P$) associated to
$X_w^P$ by the Poincar\'e duality. In what follows we will let $N
= dim(G/B)=l(w_0)$ and $N_P = dim(G/P)$.

We recall the following well known:

\begin{thm}
 The integral homology $H_*(G/P)$ is a free $\Z$--module with basis
$\{X_w^P: w \in W^P\}$.

In particular,  the integral homology $H_*(G/B)$ is a free
$\Z$-module with basis $\{X_w: w \in W \}$, where $X_w^B$ is
abbreviated by $X_w$.
\end{thm}

Since $H_*(G/P)$ is free and we have a distinguished basis (the
Schubert classes) it is reasonable to consider the basis for the
corresponding cohomology groups that are dual under the Kronecker
pairing $\<\ ,\ \>$ between homology and cohomology. Let  $\{
\epsilon_w^P : w \in W^P\}$  denote the dual basis. Thus we have
for $w,w^{\prime}\in W^P,$

$$\<\epsilon_w^P, X^P_{{w^{\prime}}}\> = \delta_{w,w^{\prime}}.$$

It suffices to study $H^*(G/B)$ because of the following well
known theorem.

\begin{thm}\label{image}
Let $\pi_P: G/B \to G/P$ be the projection. Then the  induced map
$\pi_P^*:H^*(G/P) \to H^*(G/B)$ is injective with image precisely
equal to the $W_P$-invariants of $H^*(G/B)$.

Moreover, for $w\in W^P$,
$$ \pi_P^*(\epsilon_w^P)=\epsilon_w,$$
where again we abbreviate $\epsilon_w^B$ by $\epsilon_w$.
\end{thm}

So, from now on, we will identify $H^*(G/P)$ as a subring of
$H^*(G/B)$ and denote $\epsilon_w^P$ by $\epsilon_w$ itself.

\subsection{Poincar\'e duality in $G/P$}

In the following, for each $w \in W^P$ we will need to identify
$PD(X_w^P)$ in terms of the basis $\{\epsilon_v^P\}_{v\in W^P}$.

Define the involutive map $\theta^P:W \to W$ by $\theta^P(w) =
w_0ww_{0,P}$, where $w_0$ (resp. $w_{0,P}$) is the longest element
of $W$ (resp. $W_P$).

For lack of a precise reference, we give a proof of the following:

\begin{prop}
For $w\in W^P$,  $\theta^P(w)\in W^P$  and we have
$$PD(X_w^P) = \epsilon^P_{\theta^P(w)}.$$
\end{prop}

The proposition will follow from the next two lemmas.

\begin{lem}
 The map $\theta^P$ carries
$W^P$ into itself. Moreover, we have
$$\ell(\theta^P(w)) = N_P - \ell(w),$$
where $N_P$  denotes the complex dimension of $G/P$. Thus,
$N_P=\ell(w_0)-\ell(w_{0,P})$.
\end{lem}
\proof For $w\in W^P$ and any $v\in W_P$ we have
$$\ell(w_0ww_{0,P}v)= \ell(w_0)-\ell(ww_{0,P}v).$$
But, by Lemma \ref{lengthsadd}, since $w \in W^P$ and $w_{0,P}v
\in W_P$ we get
$$\ell(ww_{0,P}v) = \ell(w) + \ell(w_{0,P}v) = \ell(w) + \ell(w_{0,P}) -
\ell(v).$$ Thus,
$$\ell(w_0ww_{0,P}v) = (\ell(w_0) - \ell(w) - \ell(w_{0,P})) + \ell(v).$$
But the above argument shows that the terms in parentheses equal
$\ell(w_0 w w_{0,P})$. Thus multiplying $\theta^P(w)$ by any $v
\in W_P$ increases the length of $\theta^P(w)$ and accordingly by
Lemma \ref{minimumlengthrepresentative} we have
$$ \theta^P(w) \in W^P.$$
Taking $v$ to be the identity in the above formula we have
$$\ell(w_0ww_{0,P}) = \ell(w_0)  - \ell(w_{0,P})- \ell(w) = N_P - \ell(w).$$
\qed

Before proving the proposition we need a general result from
algebraic topology. Let $M$ be a compact connected oriented
manifold of dimension $n$. For $a \in H_p(M)$ and $b \in
H_{n-p}(M)$ we will let $a\cdot b$ denote the intersection
pairing, \cite{Bredon},Ch.V, Section 11. We then have see
\cite{Bredon}, pg. 367,
$$\<PD(a),b\> = b \cdot a.$$

\begin{rem} In our case all the homology is even dimensional so we will
not have to worry about the interchange of order.
\end{rem}

Suppose now that $H_*(M)$ is free over $\mathbb{Z}$ and
accordingly we have $H^{p}(M) \cong \text{Hom}_\Bbb Z(H_{p}(M),
\Bbb Z)$. Thus to prove the proposition we have to identify the
element $PD(X_w^P)$ of $\text{Hom}_\Bbb Z (H_{2N_P-2\ell(w)}(G/P),
\Bbb Z)$.

The key point in the proof of the proposition is then

\begin{lem} For $v,w\in W^P$ with $\ell(v)=\ell(w)$,
$$ X_w^P \cdot X^P_{\theta^P(v)} = \delta_{w,v}.$$
\end{lem}

\proof Since the action of any element of $G$ by
left-multiplication on $G/P$ induces the trivial action on
$H^*(G/P)$, to prove the lemma it suffices to prove the following
equalities  at the cycle level:
\begin{equation}
X_w^P \cdot w_0 X^P_{\theta^P(v)} = \emptyset, \,\,\text{if}\,
v\neq w,
\end{equation}
\begin{equation}
 X_w^P \cdot w_0 X^P_{\theta^P(w)} = \{w\}.
\end{equation}
Suppose that the above cycles $X_w^P$  and $ w_0
X^P_{\theta^P(v)}$ intersect in a nonempty set. Then (since each
cycle is $T$--stable) the intersection is $T$-stable and
projective and consequently will contain a $T$-fixed point, say $u
\in W^P$. Since the Schubert cell $BwP/P$ contains the unique
$T$-fixed point $w$ and $\overline{BwP/P} = \bigsqcup\limits_{x
\leq w}BxP/P$ we find $ u\leq   w.$

Similarly, since $ u \in w_0\overline{Bw_0 v w_{0,P}P/P}$ and
$w_0$ is of order $2$, we have $w_0 uw_{0,P} \in \overline{Bw_0 v
P/P}$ and as above we find $ w_0 uw_{0,P}  \leq w_0 vw_{0,P} $ and
hence $ uw_{0,P}  \geq vw_{0,P} $. But since $u,v\in W^P$, we get
$u\geq v$ (cf. [Ku], Lemma 1.3.18).

Thus $v \leq u \leq w$. But since $\ell(v) =\ell(w)$, we obtain
$v=u=w$. By the above argument the intersection $X_w^P \cap w_0
X^P_{\theta^P(w)} = \{w\}$ set theoretically. The proof of (4) is
completed by observing that the intersection is tranverse at $w$
or alternatively by using the Poincar\'e duality and (3). \qed

Thus as operators on $H_{2N_P-2\ell(w)}(G/P)$ we have $PD(X_w^P) =
\epsilon_{\theta^P(w)}$ and the proposition follows.

\section{A formula for $\epsilon_w$ }

\subsection{The Borel model}

We continue to assume that $G$ is a connected, simply-connected
complex semisimple algebraic group. As earlier let $\omega_i$
denote the $i$-th fundamental weight. Recall that the Borel model
for $H^*(G/B)$ is obtained through the Borel homomorphism

$$\beta:\mathbb{Z}[\omega_1,\omega_2,\cdots,\omega_l] \to  H^*(G/B),  $$
which is the unique algebra homomorphism taking $\omega_i$ to the
first Chern class of the line bundle on $G/B$ associated to the
character $-\omega_i$ of $B$. It is easy to see that the
homomorphism $\beta$ commutes with the Weyl group actions and thus
for any standard parabolic subgroup $P$ of $G$, the
$W_P$-invariants
$\mathbb{Z}[\omega_1,\omega_2,\cdots,\omega_l]^{W_P}$ is mapped to
$H^*(G/P)$ under $\beta$.

Let $I \subset \mathbb{Z}[\omega_1,\omega_2,\cdots,\omega_l]$ be
the ideal generated by the $W$--invariant polynomials with zero
constant term. Then by extending the scalars to the real numbers
$\Bbb R, \beta$  induces a  surjective homomorphism (still denoted
by) $\beta:\mathbb{R}[\omega_1,\omega_2,\cdots,\omega_l] \to
H^*(G/B, \Bbb R),$ with kernel precisely equal to $I_\Bbb R:=
I\otimes_\Bbb Z \Bbb R$.

In a subsequent subsection we will use the divided-difference
operators of Demazure and Bernstein-Gelfand-Gelfand to find a
polynomial $p_w \in \mathbb{R}[\omega_1,\cdots,\omega_l]$ such
that $\beta(p_w)$ is the cohomology class $\epsilon_w$ for $w \in
W$.

\subsection{The Demazure--BGG operators}
For more details on this subsection the reader is urged to consult
\cite{Hiller}, Chapter IV and \cite{Kumar}, Chapter XI. We will
set $V= \mathfrak{a}^*$ henceforth. Let $\alpha_i$ be a simple
root and let $s_i$ be the corresponding simple reflection. We
define the divided difference operator $A_{s_i}:S^k(V) \to
S^{k-1}(V)$ by
$$A_{s_i}(f) = \frac{f - s_{i}f}{\alpha_i}.
$$

We note that $A_{s_{i}}\circ A_{s_{i}} =0$. It is also important
to note (and simple to prove)  that $A_{s_i}$ is a twisted
derivation in the following sense.

\begin{lem}
$A_{s_{i}}(pq) = A_{s_{i}}(p)q + (s_{i}p)A_{s_{i}}(q).$
\end{lem}

From the definition and the above lemma, it is easy to see that
$A_{s_{i}}$ keeps the integral form
$\mathbb{Z}[\omega_1,\omega_2,\cdots,\omega_l] \subset S(V)$
stable and, moreover, it also keeps $I_\Bbb R$ stable.

For any $w\in W$ we further define $A_w$ by
$$A_w := A_{s_1}\circ A_{s_2} \circ \cdots \circ A_{s_k}$$
where $w = s_1\cdots s_k$ is a reduced decomposition of $w$ as a
product of simple reflections. We have \cite{Hiller}, Chapter IV,
Proposition 1.7,

\begin{prop}\label{braiding}
The operators $A_w$ are well-defined, i.e., they do not depend
upon the choice of the reduced decomposition of $w$. Moreover, we
have $ A_w \circ A_v = A_{wv}$ if $\ell(wv) = \ell(w) + \ell(v)$
and $A_w \circ A_v =0$ otherwise.
\end{prop}

\subsection{The topological Demazure--BGG operators}

As earlier, for any topological space $X$, $H^*(X)$ denotes the
singular cohomology of $X$ with integral coefficients.

There is an analogue of the Demazure--BGG operator $D_{s_i}$ on
$H^*(G/B)$ (for any simple root $\alpha_i$) defined directly as
follows.

  Let $\pi_i: G/B \to G/P_i$ be the locally trivial fibration with
fibre (over $eP_i$) $P_i/B\simeq\Bbb P^1$. It is easy to see that
the restriction map $\gamma : H^*(G/B)\to H^*(P_i/B)$ is
surjective. In fact,  $\epsilon_{s_i}$ maps to the generator of
$H^2(P_i/B)$. Choose a $\Bbb Z$-module splitting $\sigma:
H^*(P_i/B)\to H^*(G/B)$ of $\gamma$. Then by the Leray-Hirsch
Theorem, the map
  $$
\Phi : H^*(G/P_i) \otimes H^*(P_i/B) \to H^*(G/B), \; u\otimes
v\mapsto (\pi^*_iu)\cup \sigma (v),
  $$
is an isomorphism.  Hence, $\pi_i^*$ is injective and $H^*(G/B)$
is a free module over $H^*(G/P_i)$ (under $\pi^*_i$) with basis 1
and $\sigma (\varepsilon)$, where $\varepsilon :=
\epsilon_{s_{i|(P_i/B)}} \in H^2(P_i/B)$ is the generator, i.e.,
  $$
H^n(G/B) \cong H^n(G/P_i) \oplus \sigma (\varepsilon )\,
H^{n-2}(G/P_i), \quad\text{ for any }n\geq 0.
  $$
Write, for any $u\in H^*(G/B)$,
  $$
u = \pi_i^*u_1 + \sigma (\varepsilon )\,\pi^*_iu_2 ,
  $$
where $u_1 \in H^*(G/P_i)$ and $u_2\in H^{*-2}(G/P_i)$ are
uniquely determined by the above equation.

Now define
  $$
D_{s_i}u := \pi_i^* u_2\in H^{*-2}(G/B).
  $$

It is easy to see that $D_{s_i}$ does not depend upon the choice
of the splitting $\sigma$.  Clearly,
  $$
D^2_{s_i} = 0.
  $$

\begin{prop} \label{intertwines}
The Borel homomorphism $\beta$ intertwines $A_{s_i}$ and $D_{s_i}$
for any simple reflection $s_i$.
\end{prop}

\proof  Let $R$ be the ring $
\mathbb{Z}[\omega_1,\omega_2,\cdots,\omega_l]$ and $R_i$  the
subring of $s_i$-invariants. Then $R$ is generated as an
$R_i$-module by $1$ and $\omega_i$. Further, $\beta (R_i) \subset
H^*(G/P_i)$ (cf. subsection 3.1). By definition, $A_{s_i}$
commutes with the multiplication by $R_i$ and also $D_{s_i}$
commutes with the multiplication by $H^*(G/P_i)$. Moreover, $
D_{s_i}(1)= A_{s_i}(1)=0.$  Thus, to prove the proposition, it
suffices to observe that
$$(*) \qquad D_{s_i} (\beta(\omega_i))=\beta(A_{s_i}(\omega_i))=1.$$
Now, it is easy to see that $\beta(\omega_i)_{\vert (P_i/B)}$ is
the generator $\varepsilon$, which proves ($*$). \qed

Thus the operators $D_{s_i}$ again satisfy the  braid relations
and we may extend $D_{s_i}$ to $D_w$ by taking a reduced
decomposition of $w$. Moreover, $D_{s_i}$ satisfies the twisted
derivation property:
$$ D_{s_i} (xy)=D_{s_i}(x)y+(s_ix)D_{s_i}(y), \,\,\text{for}\, x,y\in
H^*(G/B).
$$

We also recall the following well-known result due to Chevalley,
cf. [BGG], Theorem 3.14.
\begin{lem} \label{chevalley}
For any simple reflection $s_i$ and any $w\in W$, the cup product
\begin{equation} \epsilon_{s_i}\cdot \epsilon_w=
\sum_{w{\overset\beta\rightarrow}v} \langle w\omega_i ,
\beta^\vee\rangle\, \epsilon_v,
\end{equation}
where the notation $w{\overset\beta\rightarrow}v$ means $w\leq v$,
$\ell (v)=\ell (w)+1$, $\beta\in\Pi$ and $v=s_{\beta}w$; and
$s_{\beta}\in W$ is defined by $s_{\beta}\chi = \chi -\langle \chi
,\beta^\vee\rangle \beta $,  for $\chi\in\mathfrak{h}^*$.
 \end{lem}
The following result is of basic importance.
\begin{prop} \label{cohomologicalbgg}
For any simple reflection $s_i$ and any $w\in W$,
 $$D_{s_i} \epsilon_w=\epsilon_{ws_i}, \,\, \text{if}\,\, ws_i < w ,$$
and
$$D_{s_i} \epsilon_w=0, \,\, \text{otherwise}.$$
\end{prop}

\proof We first consider the case $ws_i > w$ so $w \in W^{P_i}$.
By Theorem \ref{image}, $\epsilon_w\in H^*(G/P_i)$, hence $D_{s_i}
\epsilon_w=0$.

So assume now that  $ws_i < w.$ By the Chevalley formula
\ref{chevalley},
$$ \epsilon_{s_i}\cdot \epsilon_{ws_i}= \epsilon_w+
\sum_{ws_i{\overset\beta\rightarrow}v, v\neq w} \langle
ws_i\omega_i , \beta^\vee\rangle\, \epsilon_v.$$ By a standard
property of Coxeter groups [Ku], Corollary 1.3.19, any $v$
appearing in the above sum satisfies $vs_i > v$. Hence applying
$D_{s_i}$ to the above equation and using the previous case, we
get
$$ D_{s_i}(\epsilon_{s_i}\cdot\epsilon_{ws_i})= D_{s_i}(\epsilon_{w}).$$
By using the twisted derivation property and the previous case
again, we get
$$ D_{s_i}(\epsilon_{s_i}\cdot\epsilon_{ws_i})= D_{s_i}(\epsilon_{s_i})
\epsilon_{ws_i}=\epsilon_{ws_i},$$ since $
D_{s_i}(\epsilon_{s_i})=1$ (because $\epsilon_{s_i}$ restricted to
$P_i/B$ equals $\varepsilon$). Combining the above two equations,
we get the proposition. \qed

\subsection{The polynomials $p_w$}

In this section we will polynomials $p_w$ such that
$$\beta(p_w) = \epsilon_w.$$

First we will find $p_{w_0}$. We define $p_{w_0}$ for $w_0$ the
longest element in $W$ as follows. Let $d$ be the product of the
positive roots. Then define

\begin{equation}\label{equationlongestelement}
p_{w_0} = \frac{d}{|W|}.
\end{equation}

\begin{prop}\label{topclass}

$$\beta(p_{w_0}) = \epsilon_{w_0}.$$
\end{prop}

As we will see that Proposition \ref{topclass} will be an almost
immediate consequence of Lemma \ref{top}. However we need a
preliminary general lemma from algebra.

Let $G$ be a finite group and $\rho:G \to$ Aut$(V)$ be a faithful
representation. Let $\mathcal{R}=S(V^*)$ be the algebra of regular
functions on $V$ and $\mathcal{F}=Q(V^*)$ be the quotient field.
The representation $\rho$ induces a representation
$\widetilde{\rho}$ from $G$ into Aut$(\mathcal{F})$ where we
consider $\mathcal{F}$ as a vector space over  the fixed field
$\mathcal{L} :=\mathcal{F}^G$. We have
$$\widetilde{\rho}(g)q(v)= q(\rho(g)^{-1}v), \,\,q(v)\in \mathcal{F}.$$

\begin{lem}\label{independence}
The set $\{ \widetilde{\rho}(g): g \in G \}$ is an independent
subset of the $\mathcal{F}$ vector space
End$_{\mathcal{L}}(\mathcal{F})$, where $\mathcal{F}$ acts on
End$_{\mathcal{L}}(\mathcal{F})$ via its multiplication on the
range.
\end{lem}
\proof Suppose we have a minimal dependence relation
$$\sum_{g \in G} q_g \widetilde{\rho}(g) = 0, q_g \in  \mathcal{F}.$$
Write $q_g = \frac{a_g}{b_g}$ with $a_g,b_g \in \mathcal{R}$ and
relatively prime. For $f \in \mathcal{R}$ we let $V(f)$ denote the
zero locus of $f$. For each $g \in G$, let $F(g)$ be the fixed
subspace of $\rho(g)$. For $g \neq I$ the subspace $F(g)$ is
proper hence $\bigcup _{g \neq I} F(g) \subsetneq V$. Choose $v_0
\in V \setminus\bigl( \bigcup_{g \neq I} F(g) \cup \bigcup_{g \in
G} V(a_g) \cup \bigcup_{g \in G} V(b_g)\bigr)$.  Then $G \to
G\cdot v_0$ is an embedding. Fix  $g_0 \in G$. Find a polynomial
$p_{g_0}$ such that
\begin{enumerate}
\item $p_{g_0}(\rho(g_0)^{-1}(v_0)) = 1$
\item $p_{g_0}(v) = 0, v \in G\cdot v_0 \setminus \{\rho(g_0)^{-1}(v_0)\}.$
\end{enumerate}
Now $0 = \sum_{g \in G} q_g(v_0) \widetilde{\rho}(g)p_{g_0}(v_0) =
\sum_{g \in G} q_g(v_0) p_{g_0}(\rho(g)^{-1}v_0) = q_{g_0}(v_0)
\neq 0$. This is a contradiction. \qed

\begin{lem}\label{top}
$A_{w_0}p_{w_0} = 1$.
\end{lem}
\proof Take a reduced decomposition $w_0 = s_{i_1}\cdots s_{i_N}$.
Then it is standard that
$\{\alpha_{i_1},s_{i_1}\alpha_{i_2},\cdots, s_{i_1}\cdots
s_{i_{N-1}}\alpha_{i_N} \}$ is an enumeration of the positive
roots. For $q \in \mathcal{F}$, let $M_q$ be the operation of
multiplication by $q$. Write $A_{w_0} = M_{\alpha_{i_1}^{-1}}(I -
s_{i_1}) \circ \cdots \circ M_{\alpha_{i_N}^{-1}}(I - s_{i_N})$.
It is then easy to see that for any $p \in \mathcal{R}$ we may
write (see the next paragraph)
\begin{equation}\label{expansion}
A_{w_0}p= \sum_{w \in W} q_w  w \cdot p.
\end{equation}
Moreover we have

\begin{equation} \label{topterm}
q_{w_0} = (-1)^N \frac{1}{d}
\end{equation}

Let $s_i$ be a simple reflection. Since $\ell(s_iw_0) < \ell(w_0)$
we have, by
 Proposition \ref{braiding}, $A_{s_i}A_{w_0} p =0$. From this we see
 that for all $p \in \mathcal{R}$ we have
 $\sum\limits _{w} q_w w\cdot p = \sum\limits_w (s_i\cdot q_w) ( s_iw \cdot
p)$.
 Apply Lemma \ref{independence} to conclude

 \begin{equation}\label{functionalequation}
 s_i\cdot q_{s_iw} = q_w.
 \end{equation}
 Combining this with  (\ref{topterm})
 we obtain
 $$q_w = (-1)^{\ell(w)+N}q_{w_0} = (-1)^{\ell(w)}\frac{1}{d}.$$
 Hence by  (\ref{expansion}) we obtain
$$ A_{w_0} p_{w_0} = \sum\limits_{w \in W} (-1)^{\ell(w)}
\frac{1}{d} w\cdot p_{w_0} =1.$$ \qed

Now we can complete the proof of Proposition \ref{topclass}. From
$dim(G/B) = N$ we deduce that   the $N$-th graded component of
$S(V^*)/I_{\Bbb R}$ is a one-dimensional vector space over $\Bbb
R$ with basis $\{\epsilon_{w_0}\}$. We will now prove that $d$ (or
$p_{w_0}$) is another basis element.

We prove that $d \notin I_{\Bbb R}$. Suppose $d \in I_{\Bbb R}$.
Then $d = \sum_i f_i g_i$ where $g_i \in S(V^*)^W$ are homogeneous
with zero constant term. By applying the alternator $\sum\limits_w
(-1)^{\ell(w)} w$ to both sides we see that we may assume that the
$f_i$ are antiinvariant elements of $S(V^*)$. But any
antiinvariant element vanishes on all the root hyperplanes and
consequently is divisible by $d$. We conclude that all the $f_i$'s
are zero. This implies that $d =0$, a contradiction.

Since $d$ has degree $N$ we find that $d$ (and hence $p_{w_0}$) is
another basis vector for the $N$-graded component of
$S(V^*)/I_{\Bbb R}$. Hence there exists $c\in \Bbb R$ such that
$p_{w_0} = c\epsilon_{w_0}$ (mod $I_\Bbb R$). But, by Proposition
\ref{cohomologicalbgg}, we have $A_{w_0} \epsilon_{w_0} =
\epsilon_1 = 1$. Hence $c= A_{w_0}p_{w_0}$. Thus, by Lemma
\ref{top}, we have $c=1$ and Proposition \ref{topclass} is proved.

We obtain as a first consequence the following

\begin{lem}\label{signaction}
The Weyl group $W$ acts on the top graded component (of degree
$\ell(w_0))$ of the graded ring $S(V)/I_\Bbb R$ by the sign
representation.
\end{lem}

We next define $p_w$ for a general $w$.
 Express $w$ as a left segment of  $w_0$.
Precisely we  find $v$ such that $w_0 = wv$ and $\ell(w_0) =
\ell(w) + \ell(v)$. Then we define
$$p_w = A_v p_{w_0} = A_{w^{-1}w_0} p_{w_0}.$$

As a consequence of Propositions \ref{intertwines}, \ref{topclass}
and \ref{cohomologicalbgg}, we get the desired realization of the
duals of the Schubert homology classes in the Borel model.

\begin{thm}\label{realization}
$$ \beta(p_w) = \epsilon_w.$$
\end{thm}

We will need the following lemma.
\begin{lem}\label{factorization}
Suppose that we have a factorization $w= w_1v_1$ with $\ell(w) =
\ell(w_1) + \ell(v_1)$.
 Then
$$A_{v_1} p_w=p_{w_1} .$$
\end{lem}
\proof Suppose we have realized $w$ as a left segment of $w_0$ by
$w_0 = wv$. Then $w_0= w_1(v_1v)$ realizes $w_1$ as a left segment
of $w_0$ with $w_1^{-1} w_0 = v_1v$ and $A_{v_1v} = A_{v_1} \circ
A_v$ by Proposition \ref{braiding}. \qed

\section{The inequalities for the rank 3 root systems}

In this section we describe the inequalities for the rank $3$ root
systems $A_3$, $C_3$ and $B_3$. Since there are many inequalities
in each case we will give only a system of representatives modulo
the action of the symmetric group $S_3$ and  leave to the reader
the task of symmetrizing the inequalities. We note that the
polytopes for $C_3$ and $B_3$ are isomorphic by Theorem
\ref{transfer}, though  we will see below that the systems are
different. We will also see that there are many trivially
redundant inequalities labelled as ($*$).

In each of the three cases there are $3$ standard maximal
parabolics, hence the system breaks up into three subsystems. We
let $r$,$s$ and $t$ be the simple reflections associated to the
nodes from left to right of the Dynkin diagram following the
Bourbaki convention (so $t$ corresponds to the long simple root in
the case of $C_3$ and the short simple root in the case of $B_3$).
In what follows $w_0$ will denote the longest element in $W$ and
$w_0^{P}$ will denote the longest element in $W^P$. Let $\lambda$
be a fundamental weight. We will  also use the notation
$X_{w\lambda}$ for the Schubert cycle $X_w^P$ with $w \in W^P$ and
$P$ the standard maximal parabolic subgroup associated to
$\lambda$.

In general the Weyl chamber  $\Delta$ is a simplicial cone, hence
in rank $3$ it is defined by $3$ linear inequalities. These
inequalities will contribute $9$ inequalities in $(v_1,v_2,v_3)$
after symmetrization.

\subsection{The inequalities for $A_3$}

In this case the quotients $G/P$ for maximal parabolics $P$ are
Grassmannians and the cohomology rings are well-known. In
particular, all the structure constants are $1$ or $0$. We will
merely record the inequalities for the three subsystems. The Weyl
chamber $\Delta$ is given by
$$\Delta = \{(x,y,z,w) : x + y + z + w = 0, x\geq y \geq z \geq w \}.$$

  We give below the inequalities  in terms of triples
$(v_1,v_2,v_3) \in \Delta^3$ with $v_i= (x_i,y_i,z_i,w_i), i
=1,2,3$.  But, to get the full set of inequalities,  we need to
symmetrize these with respect to the action of $S_3$ diagonally
permuting the variables $x_1, x_2,x_3; y_1, y_2,y_3; z_1, z_2,z_3;
w_1, w_2,w_3$.

\subsubsection{The subsystem associated to $H^*(G/P_1)$}
 In this case the quotient $G/P_1$ is $\mathbb{CP}^3$. Thus, we obtain the
subsystem (before symmetrization):
\begin{align*}
x_1 + w_2 + w_3 \leq 0 \\
y_1 + z_2 + w_3 \leq 0 \\
z_1 + z_2 + z_3 \leq 0 .\\
\end{align*}
{\it Hence there are $10$ inequalities after symmetrization.}

\subsubsection{The subsystem associated to $H^*(G/P_2)$}

In this case the quotient $G/P_2$ is the Grassmannian of $2$-
planes in $\mathbb{C}^4$. We obtain the subsystem (before
symmetrization):
\begin{align*}
x_1 + y_1 + z_2 + w_2 + z_3 + w_3 \leq 0 \\
x_1 + z_1 + y_2 + w_2 + z_3 + w_3 \leq 0 \\
x_1 + w_1 + x_2 + w_2 + z_3 + w_3 \leq 0 \\
y_1 + z_1 + y_2 + z_2 + z_3 + w_3 \leq 0 \\
\ \\
x_1 + w_1 + y_2 + w_2 + y_3 + w_3 \leq 0 \\
y_1 + z_1 + y_2 + w_2 + y_3 + w_3 \leq 0 \\
\end{align*}
{\it Hence, there are $21$ inequalities after symmetrization.}
\subsubsection{The subsystem associated to $H^*(G/P_3)$}

This subsystem is dual to the first subsystem. In this case the
quotient is the Grassmannian of $3$-planes in $\mathbb{C}^4$. We
obtain the subsystem (dual to the first)
\begin{align*}
x_1 + y_1 + z_1 + y_2 + z_2 + w_2 + y_3 + z_3 + w_3 \leq 0 \\
x_1 + y_1 + w_1 + x_2 + z_2 + w_2 + y_3 + z_3 + w_3 \leq 0 \\
\ \\
x_1 + z_1 + w_1 + x_2 + z_2 + w_2 + x_3 + z_3 + w_3 \leq 0 \\
\end{align*}
{\it Hence again there are $10$ inequalities after
symmetrization.}

{\it Thus there are altogether $50 = 41 + 9$ inequalities defining
$D_3(A_3)$ and the system is minimal, \cite{KnutsonTaoWoodward},
where $9$ in $ 41 + 9$ accounts for $9$ inequalities defining
$\Delta^3$ in $\mathfrak{a}^3$.}

\subsection{The inequalities for $C_3$}

In this subsection we take simply-connected $G$ of type $C_3$,
i.e., $G=$ Sp $(6)$.

We note that the Weyl chamber $\Delta$ is given by the triples
$x,y,z$ of real numbers satisfying
$$x\geq y \geq z \geq 0.$$
Here  $x,y,z$ are the coordinates relative to the standard basis
$\epsilon_1,\epsilon_2,\epsilon_3$ in the notation of
\cite{Bourbaki}, pg. 254 - 255. The inequalities will now be in
terms of $(v_1,v_2,v_3) \in \Delta^3$ with $v_i =
(x_i,y_i,z_i),i=1,2,3$. We will need to symmetrize the
inequalities
 with respect to the action of $S_3$ diagonally permuting the variables
$x_1, x_2,x_3; y_1, y_2,y_3; z_1, z_2,z_3$.

In what follows one will often need to verify that an expression
of an element $w \in W$ in the generators $r,s,t$ is of minimal
length. One can do this by finding the word as a connected subword
of a minimal length expression of the longest word $w_0$. Thus one
needs a plentiful supply of such expressions. The first part of
the following lemma follows from  \cite{Bourbaki}, Proposition
1.2, pg. 121. Also, in the reflection representation, the
coordinate sign changes are given by $rstsr$ (first coordinate),
$sts$ (second coordinate) and $t$ (third coordinate) and $w_0 =
-1$. From this the second part of the following lemma follows
easily.

\begin{lem}\label{longestelement}
\begin{enumerate}
\item Let $uvw$ be a product of the simple generators (in any order). Then
$(uvw)^3 = w_0$.
\item A product of the three sign-changes $rstsr,sts,t$ in any order
is equal to $w_0$.
\end{enumerate}
\end{lem}

We will also need the following proposition.

\begin{prop}
$$p_{w_0} = x^4 y^2(xyz) \,\text{mod}\,I,$$
where $I$ is the ideal generated by the $W$-invariant polynomials
in $\Bbb Z[x,y,z]$ with zero constant term.
\end{prop}
\proof By the definition of $p_{w_0}$ (cf. equation
(\ref{equationlongestelement})),
 we have
$$p_{w_0}= \frac{(x^2 - y^2)(x^2 - z^2)(y^2 - z^2)(2x)(2y)(2z)}{48}.$$
In the ring $\Z[X,Y,Z]$ we have $(X-Y)(X-Z)(Y-Z) \equiv 6 X^2Y \,
\text{mod} \, J$, where $J$ is the ideal of   $\Z[X,Y,Z]$
generated by the symmetric polynomials with zero constant term.
Now the proposition follows from the above by taking $X=x^2,
Y=y^2, Z=z^2$.\qed

We will also need the following simple fact about the Weyl group
$W$ of $C_3$. Let $w \in W$ and choose a reduced decomposition of
$w$. Let $n(w,t)$ be the number of times $t$ appears in this
decomposition.

\begin{lem}
$n(w,t)$ is independent of the reduced decomposition of $w$.
\end{lem}
\proof
From the Coxeter group property of $W$,  any one reduced
decomposition of $w\in W$ can be obtained from another by using
the Artin (i.e. generalized braid) relations, that is by replacing
$rsr$ by $srs$ or $rt$ by $tr$ or $stst$ by $tsts$. None of these
change the number of $t$'s. \qed

\subsubsection{The subsystem associated to $H^*(G/P_1)$}

We have $G/P_1 \cong \mathbb{CP}^5$ so all structure constants are
$1$ and we get one inequality for each ordered partition of $5$.
We note that $W_{P_1}$ is the group generated by $s$ and $t$.

We leave the proof of the following lemma to the reader
\begin{lem}
$W^{P_1} = \{e,r,sr,tsr,stsr,rstsr \}$.
\end{lem}
We will abbreviate the classes $\epsilon_w$ for $w \in W^{P_1}$ by
$a_i, \ i= 1,2,3,4,5$ according to the following table.  Moreover,
in the following table we list the elements $w$ of $W^{P_1}$ ,
their lengths, the {\it maximally singular weight} $\lambda_w
:=w\omega_1$ associated to the Schubert cycle
$X_{\lambda_w}=X^{P_1}_w$ and the notation for the cohomology
class $PD(X_{\lambda_w})$, the cohomology class that is Poincar\'e
dual to $X_{\lambda_w}$.

 \medskip
%\begin{figure}\label{Figure1}
\begin{center}
\begin{tabular} {|c|c|c|c|c|}
\hline
$w$    & $\ell(w)$ & $\lambda_w$ & $\epsilon_w$ & $PD(X_{\lambda_w})$ \\
\hline $e$    & $0$ & $(1,0,0)$ & $1$ & $a_5$ \\ \hline $r$    &
$1$ & $(0,1,0)$ & $a_1$ & $a_4$ \\ \hline $sr$   & $2$ & $(0,0,1)$
& $a_2$ & $a_3$ \\ \hline $tsr$  & $3$ & $(0,0,-1)$ & $a_3$ &
$a_2$ \\ \hline $stsr$ & $4$ & $(0,-1,0)$ & $a_4$ & $a_1$ \\
\hline $rstsr$ & $5$ & $(-1,0,0)$ & $a_5$ & $1$ \\ \hline

\end{tabular}
%\caption{The weights associated to $W^{P_1}$}
\end{center}
%\end{figure}

We now give the corresponding subsytem leaving to the reader the
task of symmetrizing the inequalities below. We label each
inequality with the ordered partition of $5$ that it corresponds
to. For example the label $(3,2,0)$ means the inequality
corresponds to the formula $a_3 \cdot a_2 \cdot 1 = a_5 =
\text{top class}$ in $H^*(G/P_1)$. We then refer to the above
chart  to obtain $PD(a_3)=X_{sr}^{P_1} =X_{(0,0,1)}, PD(a_2) =
X_{tsr}^{P_1} = X_{(0,0,-1)}$ and $PD(1) = X_{rstsr}^{P_1} =
X_{(-1,0,0)}$. We have
$$X_{(0,0,1)}\cdot X_{(0,0,-1)} \cdot X_{(-1,0,0)} = [pt].$$
Applying the linear functionals $(0,0,1),(0,0,-1)$ and $(-1,0,0)$
to $v_1,v_2,v_3$ respectively we get the inequality $z_1 - z_2
-x_3 \leq 0.$

 So the system of inequalities (before symmetrization) is given by:

 \begin{align*}
 x_1 \leq x_2 + x_3  \qquad (5,0,0) \\
 y_1 \leq y_2 + x_3  \qquad (4,1,0) \\
 z_1 \leq z_2 + x_3  \qquad (3,2,0) \\
 \ \\
 z_1 \leq y_2 + y_3  \qquad (3,1,1)\\
 (*) \qquad  -z_1 -z_2 -y_3 \leq 0  \qquad (2,2,1)
 \end{align*}

The  three inequalities in this subsystem  generated by ($*$)
corresponding to the ordered partition $(2,2,1)$ are trivially
redundant and do not occur in the system for $B_3$.

{\it Thus, there are $21$ inequalities in this subsystem after
symmetrization
 which includes
$3$ trivially redundant inequalities. There are no other redundant
inequalities in this subsytem.}

\subsubsection{The subsystem associated to $H^*(G/P_2)$}

The space $G/P_2$ is the space of totally-isotropic $2$-planes.
The group
 $W_{P_2}$ is  generated by the commuting simple reflections
$r$ and $t$. We leave the proof of the following lemma to the
reader.

\begin{lem}
$W^{P_2} = \{e,s,rs,ts,rts, sts,
srts,rsts,tsrts,rstrs,rtsrts,srtsrts \}$.
\end{lem}

We will abbreviate the classes $\epsilon_w$ for $w \in W^{P_2}$ to
$a_i$ or $a_i^{\prime}$ or $a_i^{\prime \prime}$ as indicated in
the next table. For the benefit of the reader, we also list for
the  elements $w$ in $W^{P_2}$, the corresponding maximally
singular weight $w\cdot \omega_2 = w\cdot (1,1,0)$ and the
Poincar\'e dual class $PD(X_w^{P_2})$.

\medskip

%\begin{figure}
 \begin{center}
\begin{tabular} {|c|c|c|c|c|}
\hline
$w$ & $\lambda_w$ & $\ell(w)$ & $\epsilon_w$  & PD($X_{\lambda_w}$) \\
\hline $e$       &   $(1,1,0)$    &  $0$ & $1$
& $a_7$
\\ \hline
$s$       &   $(1,0,1)$    &  $1$ & $a_1$                   &
$a_6$
\\ \hline
$rs$      &   $(0,1,1)$    &  $2$ & $a_2^{\prime}$          &
$a_5^{\prime}$
\\ \hline
$ts$      &   $(1,0,-1)$   &  $2$ & $a_2^{\prime \prime}$   &
$a_5^{\prime \prime}$    \\ \hline $rts$     &   $(0,1,-1)$   &
$3$ & $a_3^{\prime}$          & $a_4^{\prime}$
\\ \hline
$sts$     &   $(1,-1,0)$   &  $3$ & $a_3^{\prime \prime}$    &
$a_4^{\prime \prime}$    \\ \hline $srts$    &   $(0,-1,1)$   &
$4$ & $a_4^{\prime}$          & $a_3^{\prime}$
\\ \hline
$rsts$    &   $(-1,1,0)$   &  $4$ & $a_4^{\prime \prime}$    &
$a_3^{\prime \prime}$    \\ \hline $tsrts$   &   $(0,-1,-1)$  &
$5$ & $a_5^{\prime}$          & $a_2^{\prime}$
\\ \hline
$rstrs$   &   $(-1,0,1)$   &  $5$ & $a_5^{\prime \prime}$    &
$a_2^{\prime \prime}$    \\ \hline $rtsrts$  &   $(-1,0,-1)$  &
$6$ & $a_6$                   &  $a_1$
\\ \hline
$srtsrts$ &   $(-1,-1,0)$  &  $7$ & $a_7$                   &  $1$
\\ \hline
\end{tabular}
%\caption{The weights associated to $W^{P_2}$}
\end{center}
%\end{figure}

 We list the polynomials $p_w$
(mod $I$) in the next table, which is obtained by applying Lemma
\ref{factorization} and observing that both of $x^4y^2+ x^2y^4$
and $x^4+y^4+x^2y^2$ belong to $I$.

In the following table it will be convenient to use $g_n$ to
denote the finite geometric series
$$ g_n = \sum_{i=0}^n x^i y ^{n-i}.$$

\begin{center}
\begin{tabular} {|c|c|c|c|}
\hline $w$    &  $v$  & $\tilde{w}=wv$ & $p_w (\text{mod}\,I)$ \\
\hline $srtsrts$    & $rt$  & $w_0$ & $(xy)^3 (x + y)$
\\ \hline $rtsrts$    & $rts$ & $w_0$  &  $(xy)^3$
\\ \hline $tsrts$   & $rtsr$ & $w_0$   &  $x y (x +y)(x y + z^2)$
\\ \hline $rstrs$  & $trst$  & $w_0$   &  $x y g_3$
\\ \hline $srts$ & $rtsrt$   & $w_0$   &  $(xy)^2$
\\ \hline $rsts$ & $rstst$  & $w_0$    &  $g_4$
\\ \hline $rts$ &  $rts$    & $rtsrts$ &  $x y (x +y)$
\\ \hline $sts$ & $trstsr$   & $w_0$   &  $g_3$
\\ \hline $rs$ & $ts$       & $rsts$   &  $xy$
\\ \hline
$ts$ & $trs$      & $ tstrs$ &  $g_2$                               \\
\hline $s$ & $ts$        & $ sts$   &  $x + y$
\\ \hline

\end{tabular}
\end{center}

Now computing the products of $p_w$ (mod $I$) and using the
Chevalley formula Lemma \ref{chevalley}, we obtain the following.
\begin{thm}

The cohomology ring $H^*(G/P_2)$ is given by the following table.
\begin{center}
\begin{tabular} {|c|c|c|c|c|c|}
\hline $H^*(G/P_2)$ & $a_1$ & $a_2^{\prime}$ & $a_2^{\prime
\prime}$ & $a_3^{\prime}$ & $a_3^{\prime \prime}$  \\ \hline

$a_1$ & $a_2^{\prime} + a_2^{\prime \prime}$ & $a_3^{\prime}$ &
$a_3^{\prime} + a_3^{\prime \prime}$ & $ 2 a_4^{\prime} +
a_4^{\prime \prime}$ & $ a_4^{\prime} + 2 a_4^{\prime \prime}$ \\
\hline

$a_2^{\prime}$ & & $a_4^{\prime}$ & $a_4^{\prime} + a_4^{\prime
\prime}$ & $ a_5^{\prime} + a_5^{\prime \prime}$ & $ a_5^{\prime
\prime}$  \\ \hline

$a_2^{\prime \prime}$ & & & $ 2 a_4^{\prime} + 2 a_4^{\prime
\prime}$ & $a_5^{\prime} + 2 a_5^{\prime \prime}$ & $ a_5^{\prime}
+ 2 a_5^{\prime \prime} $  \\ \hline

$a_3^{\prime }$ & & & & $ 2 a_6$ & $a_6$ \\ \hline

$a_3^{\prime \prime}$ & & & & & $ 2 a_6$ \\ \hline

\end{tabular}
\end{center}

\begin{center}
\begin{tabular} {|c|c|c|c|c|c|c|}
\hline $H^*(G/P_2)$ & $a_4^{\prime}$ & $a_4^{\prime \prime}$ &
$a_5^{\prime}$ & $a_5^{\prime \prime}$ & $a_6$ & $a_7$ \\ \hline

$a_1$ & $a_5^{\prime}+a_5^{\prime \prime}$ & $a_5^{\prime \prime}$
& $a_6$ & $a_6$ & $a_7$ & $0$ \\ \hline

$a_2^{\prime}$ & $a_6$ & $0$ & $a_7$ & $0$ & $0$ & $0$ \\ \hline

$a_2^{\prime \prime}$ & $a_6$ & $a_6$ & $0$ & $a_7$ & $0$ & $0$ \\
\hline

$a_3^{\prime }$ & $a_7$ & $0$ & $0$ & $0$ & $0$ & $0$\\ \hline

$a_3^{\prime \prime}$ &  $0$ & $a_7$ & $0$ & $0$ & $0$ & $0$
\\ \hline
\end{tabular}
\end{center}
\end{thm}

We now write down the subsystem (before symmetrization).

\begin{align*}
x_1 + y_1 & \leq x_2 + y_2 + x_3 + y_3  \qquad    & (7,0,0) \\
x_1 + z_1 & \leq x_2 + z_2 + x_3 + y_3  \qquad    & (6,1,0) \\
y_1 + z_1 & \leq y_2 + z_2 + x_3 + y_3  \qquad    &
(5^{\prime},2^{\prime},0) \\
x_1 + z_2 & \leq z_1 + x_2 + x_3 + y_3  \qquad    & (5^{\prime
\prime},
2^{\prime \prime},0)\\
y_1 + z_2 & \leq z_1 + y_2 + x_3 + y_3  \qquad    &
(4^{\prime},3^{\prime},0) \\
x_1 + y_2 & \leq y_1 + x_2 + x_3 + y_3  \qquad    & (4^{\prime
\prime},
3^{\prime \prime},0) \\
\ \\
y_1 + z_1 & \leq x_2 + z_2 + x_3 + z_3  \qquad    & (5^{\prime},1,1) \\
(**)\qquad x_1 \leq z_1 & + x_2 + z_2 + x_3 + z_3   \qquad   &
(5^{\prime
\prime}, 1,1) \\
(**)\qquad y_1 \leq z_1 & + y_2 + z_2 + x_3 + z_3   \qquad   &
(4^{\prime},2^{\prime},1) \\
y_1 + z_2 & \leq z_1 + x_2 + x_3 + z_3   \qquad   &
(4^{\prime},2^{\prime
\prime},1) \\
x_1 + z_2 & \leq y_1 + x_2 + x_3 + z_3   \qquad   & (4^{\prime
\prime},2^{\prime \prime},1,) \\
(*)\qquad y_1 + z_2 & \leq x_1 + y_2 + x_3 + z_3   \qquad   &
(3^{\prime
\prime},3^{\prime},1) \\
\ \\
(*)\qquad z_1 \leq y_1 & + y_2 + z_2 + y_3 + z_3   \qquad   &
(3^{\prime},2^{\prime},2^{\prime}) \\
(*)\qquad z_1 + z_2 & \leq y_1 + x_2 + y_3 + z_3   \qquad   &
(3^{\prime},2^{\prime \prime},2^{\prime}) \\
(*)\qquad y_1 + z_2 & \leq x_1 + x_2 + y_3 + z_3   \qquad   &
(3^{\prime
\prime},2^{\prime \prime},2^{\prime}) \\
\end{align*}

We explain how an inequality corresponds to a decorated ordered
partition by the example of $(4^{\prime},2^{\prime \prime},1)$.
The decorated partition corresponds to the three-fold product
$a_4^{\prime}\cdot a_2^{\prime \prime}\cdot a_1 = a_7 = \text{top
class}$ in $H^*(G/P_2)$. Taking the Poincar\'e dual cycles, the
above product corresponds to the intersection product
$X_{rts}^{P_2}\cdot X_{rstrs}^{P_2}\cdot X_{rtsrts}^{P_2} = [pt]$.
Using the correspondence between Weyl group elements in $W^{P_2}$
and linear functionals arising from  the Weyl group orbit of
$(1,1,0)$, we find that the three Weyl group elements indexing the
cycles in the intersection product correspond to the linear
functionals $(0,1,-1),(-1,0,1),(-1,0,-1)$. Applying these linear
functionals to $v_1,v_2,v_3$ respectively and collecting terms one
obtains the inequality $y_1 -z_1 - x_2 + z_2 -x_3 -z_3 \leq 0$ or
equivalently $y_1 + z_2 \leq z_1 + x_2 + x_3 + z_3$.

The inequalities in the subsystem labelled ($*$) corresponding to
the ordered partitions $(3^{\prime
\prime},3^{\prime},1),(3^{\prime},2^{\prime},2^{\prime}),
(3^{\prime},2^{\prime \prime},2^{\prime}), (3^{\prime
\prime},2^{\prime \prime},2^{\prime}) $ are trivially redundant.
They generate (after symmetrization) respectively $6$,$3$,$6$,$6$
trivially redundant inequalities.

The inequalities corresponding to the ordered partitions
$(5^{\prime \prime},1,1)$, $(4^{\prime},2^{\prime},1)$ and
$(3^{\prime}, 2^{\prime},2^{\prime})$ do not occur in the system
for $B_3$. Consequently they too must be redundant. We now check
this directly.

First the three inequalities corresponding to the decorated
ordered partition $(3^{\prime},2^{\prime}, 2^{\prime})$ are
trivially redundant. In order to check that the three inequalities
corresponding to $(5^{\prime \prime},1,1)$ are redundant, we
observe that we have $x_1\leq x_2 + x_3$ from the first subsystem
(corresponding to $G/P_1$). As a consequence we have $x_1 \leq x_2
+ x_3 + z_1 + z_2 + z_3$. Finally to check that the six
inequalities corresponding to $(4^{\prime}, 2^{\prime},1)$ are
redundant, we observe that we have $y_1 \leq y_2 + x_3$ from the
first subsystem again. Hence $y_1 \leq y_2 + x_3 + z_1 + z_2 +
z_3$.

{\it This subsystem (corresponding to $G/P_2$) after
symmetrization consists of $78$ inequalities of which $21$ are
trivially redundant (marked by ($*$)).  There are $9$ more
redundant inequalities marked by ($**$). These $9$ inequalities do
not occur in the system for $B_3$.}

\subsubsection{The subsystem associated to $H^*(G/P_3)$}

The space $G/P_3$ is the space of totally-isotropic $3$-planes,
i.e., the Lagrangian Grassmannian. The group  $W_{P_3}$ is
generated by the simple reflections $r$ and $s$. We leave the
proof of the following lemma to the reader.

\begin{lem}
$$W^{P_3}= \{e,t,st,rst,tst,trst,strst,tstrst\}.$$
\end{lem}

We again abbreviate the classes $\epsilon_w$ to $a_i$ or
$a_i^{\prime}$ or $a_i^{\prime \prime}$ according to the following
table.

\medskip
%\begin{figure}
 \begin{center}
\begin{tabular} {|c|c|c|c|c|}
\hline
$w$       & $\ell(w)$ & $\lambda_w$ & $\epsilon_w$ & $PD(X_{\lambda_w})$ \\
\hline
 $e$       & $0$ &  $(1,1,1) $    & $1$           & $a_6$                \\
\hline
 $t$       & $1$ &  $(1,1,-1)$    & $a_1$               & $a_5$
\\ \hline
 $st$      & $2$ &  $(1,-1,1)$    & $a_2$            & $a_4$
\\ \hline
 $rst$     & $3$ &  $(-1,1,1)$    & $a_3^{\prime}$             &
$a_3^{\prime \prime}$          \\ \hline
 $tst$     & $3$ &  $(1,-1,-1)$   & $a_3^{\prime \prime}$         &
$a_3^{\prime}$    \\ \hline
 $trst$    & $4$ &  $(-1,1,-1)$   & $a_4 $          & $a_2$
\\ \hline
 $strst$   & $5$ &  $(-1,-1,1)$   & $a_5$         & $a_1 $                \\
\hline
 $tstrst$  & $6$ &  $(-1,-1,-1)$  & $a_6$          & $ 1$
\\ \hline

\end{tabular}
%\caption{The weights associated to  $W^{P_3}$}
\end{center}
%\end{figure}

We again give the formulas for $p_w$ (mod $I$). In the following
table we define a symmetric cubic polynomial $f(x,y,z)$ by
$f(x,y,z) = x^2y + xy^2 + x^2z + xz^2 + y^2 z + yz^2$.

\begin{center}
\begin{tabular} {|c|c|c|c|}
\hline $w$    &  $v$  & $\tilde{w}=wv$ & $p_w (\text{mod}\,I)$ \\
\hline $w_0^{P_3}=tstrst$    & $rsr$  & $w_0$ & $(x y z)f $
\\ \hline $strst$    & $rstr$ & $w_0$  &  $(x y z)(xy + xz + yz)$
\\ \hline
$trst$   & $rstrs$ & $w_0$   &  $(x y z) (x +y + z)$  \\ \hline
$tst$  & $rst$  & $w_0^{P_3}= tstrst$   &  $f + x y z$
\\ \hline
$rst$ & $rstrst$   & $w_0$   &  $x y z$                       \\
\hline $st$ & $rst$   & $strst$   &  $xy + x z + y z$         \\
\hline $t$ & $st$        & $ tst$   &  $x + y + z$      \\ \hline

\end{tabular}
\end{center}

We now have the following  table  constructed by multiplying the
polynomials $p_w$ modulo $I$.

\begin{thm}

The cohomology ring $H^*(G/P_3)$ is given by the following table.

\begin{center}
\begin{tabular} {|c|c|c|c|c|c|c|c|}
\hline

$H^*(G/P_3)$ & $a_1$ &  $a_2$ &  $a_3^{\prime}$ & $a_3^{\prime
\prime}$ & $a_4$ & $a_5$ &  $a_6$   \\ \hline $a_1$ & $2 a_2$ & $2
a_3^{\prime} + a_3^{\prime \prime}$ & $a_4$ & $2a_4$ & $2a_5$ &
$a_6$ & $0$ \\ \hline $a_2$ & & $2a_4$ & $a_5$ & $2a_5$ & $a_6$ &
$0$ & $0$ \\ \hline $a_3^{\prime}$ & & & $0$ & $a_6$ & $0$ & $0$ &
$0$ \\ \hline $a_3^{\prime \prime}$ & & & & $0$ & $0$ & $0$ & $0$
\\ \hline

\end{tabular}
\end{center}
\end{thm}

Using the above tables, the reader can easily verify that we have
the following inequalities (which have to be symmetrized).

\begin{align*}
x_1 + y_1 + z_1 \leq x_2 + y_2 + z_2 + x_3 + y_3 + z_3 \qquad & (6,0,0) \\
x_1 + y_1 + z_2 \leq z_1 + x_2 + y_2 + x_3 + y_3 + z_3 \qquad & (5,1,0) \\
x_1 + z_1 + y_2 \leq y_1 + x_2 + z_2 + x_3 + y_3 + z_3 \qquad & (4,2,0) \\
x_1 + y_2 + z_2 \leq y_1 + z_1 + x_2 + x_3 + y_3 + z_3 \qquad &
(3^{\prime},3^{\prime \prime},0)\\
\ \\
x_1 + y_2 + z_3 \leq y_1 + z_1 + x_2 + z_2 + x_3 + y_3 \qquad &
(3^{\prime},2,1) \\
\end{align*}

{\it This gives that the subsystem corresponding to $G/P_3$
consists of $27$ inequalities. None of them are trivially
redundant.}

 The $27$ inequalities above can be rewritten  in a very simple way.
 Let $S = \sum_{i=1}^3 x_i + y_i +z_i$. Then the $27$ inequalities
 are just the inequalities
 $$x_i + y_j + z_k \leq \frac{S}{2}, \  i,j,k = 1,2,3.$$

{\it Thus finally we see that for $C_3$ there are $135 = 126 + 9$
inequalities of which $24$ are trivially redundant  ($9$ in
$126+9$ coming from the inequalities defining $\Delta^3$ inside
$\mathfrak{a}^3$). There are $9$ more redundant inequalities.
These $9$ inequalities do not occur in the system for $B_3$ .
Hence the  subsystem for $C_3$ can be brought down to  altogether
$102$ inequalities. Moreover, a computer calculation shows that
the polyhedral cone $D_3(C_3)$ has exactly $102$ faces and thus
these $102$ inequalities are irredundant.}

\subsection{The inequalities for $B_3$}

In this subsection we take simply-connected $G$ of type $B_3$,
i.e., $G=$ Spin $(7)$.

We note that the Weyl chamber $\Delta$ is given by triples $x,y,z$
of real numbers satisfying
$$x\geq y \geq z \geq 0.$$

The inequalities will now be in terms of $(v_1,v_2,v_3) \in
\Delta^3$ with $v_i = (x_i,y_i,z_i),i=1,2,3$.

The Weyl groups for $\text{Sp}(6)$ and $\text{Spin}(7)$ are
isomorphic. In fact, in the standard coordinates $(x,y,z)$ they
are identical. Hence the sets $W^{P_i},i = 1,2,3$ will be
identical. However the operators $A_w$ and the polynomials $p_w$
will be different (but proportional) as we now see.
 To distinguish, we will denote them by $p_w^{\text{Sp} (6)}$ and
$p_w^{\text{Spin} (7)}$ respectively.

We will need the following proposition.

\begin{prop}
$$p_{w_0}^{\text{\rm Spin} (7)} = \frac{x^4 y^2(xyz)}{8}\,\,(\text{\rm
mod}\,I).$$
\end{prop}

\proof By equation (\ref{equationlongestelement}) we have
$$p_{w_0}^{\text{Spin} (7)}= \frac{xyz (x^2 - y^2)(x^2 - z^2)(y^2 - z^2)
}{48}.$$ Hence the proposition follows from the corresponding
result for $\text{Sp}(6)$. \qed

The next lemma tells us how to read off the polynomials $p_w$ for
the case of $\text{Spin}(7)$ from the corresponding polynomials
for $\text{Sp}(6)$. We will temporarily use the notation
$A_w^{\text{Spin}(7)}$ and $A_w^{\text{Sp}(6)}$ for the
Demazure-BGG operators $A_w$ for the groups Spin$(7)$ and Sp$(6)$
respectively.

\begin{lem}

Let $v\in W$ and let $n(v,t)$ be the number of times the simple
reflection $t$ occurs in some reduced decomposition of $v$. Then
we have
$$A_v^{\text{\rm Spin}(7)} = 2^{n(v,t)} A_v^{\rm Sp(6)}.$$

\end{lem}

\proof $A_t^{\text{Spin}(7)}=2 \ A_t^{Sp(6)}$ and $
A_r^{\text{Spin}(7)} = A_r^{Sp(6)}$,  $A_s^{\text{Spin}(7)}=
A_s^{Sp(6)}.$ \qed

\begin{cor}\label{proportionality} For any $w\in W$,
$p_w^{\text{\rm Spin}(7)}=\frac{2^{n(w^{-1}w_0,t)}}{8} \  p_w^{\rm
Sp(6)} = 2^{- n(w,t)}  \  p_w^{\rm Sp(6)}$.
\end{cor}
\proof Note that $n(w,t) = 8 - n(v,t)$, for $v=w^{-1}w_0$. \qed

\subsubsection{The subsystem associated to $H^*(G/P_1)$}

We note that $W_{P_1}$ is the group generated by $s$ and $t$. We
have (since this is what we had for $\text{Sp}(6)$)

\begin{lem}
$W^{P_1}=\{e,r,sr,tsr,stsr,rstsr \}.$
\end{lem}

We have $G/P_1 \cong Q_5$, the smooth quadric hypersurface in
$\mathbb{CP}^6$   so the inequalities will be parametrized by a
subset of the  ordered partitions of $5$. We will abbreviate the
classes $\epsilon_w$ to $a_i$ according to the following table. In
addition, we list the elements $w$ of $W^{P_1}$ , their lengths,
the maximally singular weight $\lambda_w$ associated to the
Schubert cycle $X_w^{P_1}=X_{\lambda_w}$ and the notation for the
cohomology class $PD(X_{\lambda_w})$, the cohomology class that is
Poincar\'e dual to $X_{\lambda_w}$.

\begin{center}
\begin{tabular} {|c|c|c|c|c|}
\hline
$w$    & $\ell(w)$ & $\lambda_w$ &$\epsilon_w$ & $PD(X_{\lambda_w})$ \\
\hline $e$    & $0$ & $(1,0,0)$ & $1$ & $a_5$ \\ \hline $r$    &
$1$ & $(0,1,0)$ & $a_1$      & $a_4$ \\ \hline $sr$   & $2$ &
$(0,0,1)$  & $a_2$      & $a_3$ \\ \hline $tsr$  & $3$ &
$(0,0,-1)$  & $a_3$             & $a_2$ \\ \hline $stsr$ & $4$ &
$(0,-1,0)$  & $a_4$             & $a_1$ \\ \hline $rstsr$ & $5$ &
$(-1,0,0)$ & $a_5$              & $1$ \\ \hline

\end{tabular}
\end{center}

The following theorem follows easily from Corollary
\ref{proportionality} by using the corresponding theorem for
Sp$(6)$.
\begin{thm}

The multiplication table for $H^*(G/P_1)$ is given by the
following table.

\begin{center}
\begin{tabular} {|c|c|c|c|c|c|}
\hline

$H^*(G/P_1)$ & $a_1$ &  $a_2$ &  $a_3$ &  $a_4$ & $a_5$   \\
\hline $a_1$ & $ a_2$ & $2 a_3$ & $a_4$ &  $a_5$ &  $0$ \\ \hline
$a_2$ & & $2a_4$ & $a_5$ &  $0$ & $0$ \\ \hline

\end{tabular}
\end{center}
\end{thm}

We now give the corresponding subsystem, leaving to the reader the
task of symmetrizing the inequalities below.

 \begin{align*}
 x_1 \leq x_2 + x_3 \qquad & (5,0,0) \\
 y_1 \leq y_2 + x_3  \qquad & (4,1,0) \\
 z_1 \leq z_2 + x_3 \qquad & (3,2,0) \\
 z_1 \leq y_2 + y_3  \qquad & (3,1,1) \\
 \end{align*}

{\it After symmetrizing there are $18$ inequalities, none are
trivially redundant.}

\subsubsection{The subsystem associated to $H^*(G/P_2)$}

The space $G/P_2$ is the space of totally-isotropic $2$-planes. We
have that $W_{P_2}$ is the group generated by the commuting simple
reflections $r$ and $t$. We have, as for $\text{Sp}(6)$,

\begin{lem}
$W^{P_2} = \{e,s,rs,ts,rts, sts,
srts,rsts,tsrts,rstrs,rtsrts,srtsrts \}$.
\end{lem}

We will abbreviate the classes $\epsilon_w$ by $b_i$ or
$b_i^{\prime}$ or $b_i^{\prime \prime}$ according to the following
table.

The following table follows easily from the corresponding tables
for Sp$(6)$ and Corollary \ref{proportionality}.

 \begin{center}
\begin{tabular} {|c|c|c|c|c|c|}
\hline $w$ & $\lambda_w$ & $\ell(w)$  & $\epsilon_w$ & $p_w
\,(\text{mod}\, I)$ & $PD(X_{\lambda_w}$) \\ \hline $e$       &
$(1,1,0)$    &  $0$ & $1$                & $1$ & $b_7$
\\ \hline $s$       &   $(1,0,1)$    &  $1$ & $b_1$
&  $x + y$                              & $b_6$
\\ \hline $rs$      &   $(0,1,1)$    &  $2$ & $b_2^{\prime}$
&  $ x y $                               &  $b_5^{\prime}$
\\ \hline $ts$      &   $(1,0,-1)$   &  $2$ & $b_2^{\prime
\prime}$              &
$1/2 \  g_2$                                 &  $b_5^{\prime \prime}$    \\
\hline $rts$     &   $(0,1,-1)$   &  $3$ & $b_3^{\prime}$
& $1/2 \  x y (x +y)$                  &  $b_4^{\prime}$
\\ \hline $sts$     &   $(1,-1,0)$   &  $3$ & $b_3^{\prime
\prime}$               &
$1/2 \  g_3$                                 &  $b_4^{\prime \prime}$    \\
\hline $srts$    &   $(0,-1,1)$   &  $4$ & $b_4^{\prime}$
&
$1/2 \  (x y)^2$                         &  $b_3^{\prime}$           \\
\hline $rsts$    &   $(-1,1,0)$   &  $4$ & $b_4^{\prime \prime}$
&
$1/2 \  g_4$                                 &  $b_3^{\prime \prime}$    \\
\hline $tsrts$   &   $(0,-1,-1)$  &  $5$ & $b_5^{\prime}$
& $1/4 \ x y (x +y)(x y + z^2)$ &  $b_2^{\prime}$           \\
\hline $rstrs$   &   $(-1,0,1)$   &  $5$ & $b_5^{\prime \prime}$
&
$1/2 \ x y g_3$                         &  $b_2^{\prime \prime}$    \\
\hline $rtsrts$  &   $(-1,0,-1)$  &  $6$ & $b_6$
&
$1/4 \ (xy)^3$                          &  $b_1$                    \\
\hline $srtsrts$ &   $(-1,-1,0)$  &  $7$ & $b_7$
& $1/4 \ (xy)^3 (x + y)$              &  $1$
\\ \hline
\end{tabular}
\end{center}

From these formulae one can easily construct the multiplication
table for $H^*(G/P_2)$.

\begin{thm}
The cohomology ring $H^*(G/P_2)$ is given by the following table.

\begin{center}
\begin{tabular} {|c|c|c|c|c|c|}
\hline $H^*(G/P_2)$ & $b_1$ & $b_2^{\prime}$ & $b_2^{\prime
\prime}$ & $b_3^{\prime}$ & $b_3^{\prime \prime}$ \\ \hline

$b_1$ & $b_2^{\prime} + 2 b_2^{\prime \prime}$ & $2 b_3^{\prime}$
& $b_3^{\prime} + b_3^{\prime \prime}$ & $ 2 b_4^{\prime} +
b_4^{\prime \prime}$ & $ b_4^{\prime} + 2 b_4^{\prime \prime}$  \\
\hline

$b_2^{\prime}$ & & $2 b_4^{\prime}$ & $b_4^{\prime} + b_4^{\prime
\prime}$ & $2 b_5^{\prime} + b_5^{\prime \prime}$
& $ b_5^{\prime \prime}$ \\
\hline

$b_2^{\prime \prime}$ & & & $ b_4^{\prime} + b_4^{\prime \prime}$
& $b_5^{\prime} + b_5^{\prime \prime}$ & $ b_5^{\prime} +
b_5^{\prime \prime} $  \\ \hline

$b_3^{\prime }$ & & & & $ 2 b_6$ & $b_6$ \\ \hline

$b_3^{\prime \prime}$ & & & & & $ 2 b_6$ \\ \hline

\end{tabular}
\end{center}

\begin{center}
\begin{tabular} {|c|c|c|c|c|c|c|}
\hline $H^*(G/P_2)$ & $b_4^{\prime}$ & $b_4^{\prime \prime}$ &
$b_5^{\prime}$ & $b_5^{\prime \prime}$ & $b_6$ & $b_7$ \\ \hline

$b_1$ & $2 b_5^{\prime}+b_5^{\prime \prime}$ & $b_5^{\prime
\prime}$ & $b_6$ & $2 b_6$ & $b_7$ & $0$ \\ \hline

$b_2^{\prime}$ & $2 b_6$ & $0$ & $b_7$ & $0$ & $0$ & $0$ \\
\hline

$b_2^{\prime \prime}$ & $b_6$ & $b_6$ & $0$ & $b_7$ & $0$ & $0$ \\
\hline

$b_3^{\prime }$ & $b_7$ & $0$ & $0$ & $0$ & $0$ & $0$\\ \hline

$b_3^{\prime \prime}$ &$0$ & $b_7$ & $0$ & $0$ & $0$ & $0$
\\ \hline

\end{tabular}
\end{center}

\end{thm}

We then read off the inequalities.

\begin{align*}
x_1 + y_1 & \leq x_2 + y_2 + x_3 + y_3 \qquad & (7,0,0) \\
x_1 + z_1 & \leq x_2 + z_2 + x_3 + y_3 \qquad & (6,1,0) \\
y_1 + z_1 & \leq y_2 + z_2 + x_3 + y_3 \qquad & (5^{\prime},2^{\prime},0) \\
x_1 + z_2 & \leq z_1 + x_2 + x_3 + y_3 \qquad & (5^{\prime
\prime},2^{\prime
\prime},0)\\
y_1 + z_2 & \leq z_1 + y_2 + x_3 + y_3 \qquad & (4^{\prime},3^{\prime},0)\\
x_1 + y_2 & \leq y_1 + x_2 + x_3 + y_3 \qquad & (4^{\prime
\prime},
3^{\prime \prime},0)\\
\ \\
y_1 + z_1 & \leq x_2 + z_2 + x_3 + z_3 \qquad & (5^{\prime},1,1) \\
y_1 + z_2 & \leq z_1 + x_2 + x_3 + z_3 \qquad &
(4^{\prime},2^{\prime
\prime},1,) \\
x_1 + z_2 & \leq y_1 + x_2 + x_3 + z_3 \qquad & (4^{\prime \prime}
2^{\prime
\prime},1)\\
(*)\qquad y_1 + z_2 & \leq x_1 + y_2 + x_3 + z_3 \qquad &
(3^{\prime
\prime},3^{\prime},1) \\
\ \\
(*)\qquad y_1 + z_2 & \leq x_1 + x_2 + y_3 + z_3 \qquad &
(3^{\prime
\prime},2^{\prime \prime},2^{\prime})\\
(*)\qquad z_1 + z_2 & \leq y_1 + x_2 + y_3 + z_3 \qquad &
(3^{\prime},2^{\prime \prime},2^{\prime}) \\
\ \\
(*)\qquad z_1 + z_2 & + z_3 \leq y_1 + x_2 + x_3 \qquad &
(3^{\prime},2^{\prime \prime},2^{\prime \prime})\\
(*)\qquad y_1 + z_2 & + z_3 \leq x_1 + x_2 + x_3 \qquad &
(3^{\prime
\prime},2^{\prime \prime},2^{\prime \prime})\\
\end{align*}

The inequalities corresponding to the ordered partitions
$(3^{\prime \prime},3^{\prime}, 1)$ $(3^{\prime \prime},2^{\prime
\prime}, 2^{\prime})$, $(3^{\prime },2^{\prime
\prime},2^{\prime})$, $(3^{\prime },2^{\prime \prime},2^{\prime
\prime})$ and  $(3^{\prime \prime},2^{\prime \prime},2^{\prime
\prime})$ are trivially redundant. After symmetrizing they give
rise to $24$ trivially redundant inequalities. The trivially
redundant inequalities corresponding to the decorated ordered
partitions $(3^{\prime},2^{\prime \prime}, 2^{\prime \prime})$ and
$(3^{\prime \prime},2^{\prime \prime},2^{\prime \prime})$ do not
occur in the system for $\text{Sp}(6)$.

{\it There are $72$ inequalities after symmetrizing, of which
 $24$ are trivially redundant.}

\subsubsection{The subsystem associated to $H^*(G/P_3)$}

The space $G/P_3$ is the space of totally-isotropic $3$-planes in
$\mathbb{C}^7$. The group  $W_{P_3}$ is generated by the simple
reflections $r$ and $s$. We have the following from the
corresponding result for Sp$(6)$.

\begin{lem}
$$W^{P_3}= \{e,t,st,rst,tst,trst,strst,tstrst\}.$$

\end{lem}

We again abbreviate the classes $\epsilon_w$ to $b_i$'s or
$b_i^{\prime}$'s or $b_i^{\prime \prime}$'s  according to the
following table.

 \begin{center}
\begin{tabular} {|c|c|c|c|c|}
\hline $w$       & $\ell(w)$ & $\lambda_w$ & $\epsilon_w$  &
$PD(X_{\lambda_w})$
\\ \hline
 $e$       & $0$ &  $(1,1,1) $     & $1$   &$b_6$            \\ \hline
 $t$       & $1$ &  $(1,1,-1)$     & $b_1$  & $b_5$               \\ \hline
 $st$      & $2$ &  $(1,-1,1)$     & $b_2$ & $b_4$                \\ \hline
 $rst$     & $3$ &  $(-1,1,1)$     & $b_3^{\prime}$ & $b_3^{\prime \prime}$
\\ \hline
 $tst$     & $3$ &  $(1,-1,-1)$    & $b_3^{\prime \prime}$  & $b_3^{\prime}$
\\ \hline
 $trst$    & $4$ &  $(-1,1,-1)$    & $b_4$     & $b_2$                 \\
\hline
 $strst$   & $5$ &  $(-1,-1,1)$    & $b_5$   & $b_1 $                \\
\hline
 $tstrst$  & $6$ &  $(-1,-1,-1)$   & $b_6$   & $ 1$                   \\
\hline

\end{tabular}
\end{center}

\begin{thm}

The cohomology ring is given by the following table.

\begin{center}
\begin{tabular} {|c|c|c|c|c|c|c|c|}
\hline

$H^*(G/P_3)$ & $b_1$ &  $b_2$ &  $b_3^{\prime}$ & $b_3^{\prime
\prime}$ & $b_4$ & $b_5$ &  $b_6$   \\ \hline $b_1$ & $ b_2$ & $
b_3^{\prime} + b_3^{\prime \prime}$ & $b_4$ & $b_4$ & $b_5$ &
$b_6$ & $0$ \\ \hline $b_2$ & & $2b_4$ & $b_5$ & $b_5$ & $b_6$ &
$0$ & $0$ \\ \hline $b_3^{\prime}$ & & & $0$ & $b_6$ & $0$ & $0$ &
$0$ \\ \hline $b_3^{\prime \prime}$ & & & & $0$ & $0$ & $0$ & $0$
\\ \hline

\end{tabular}
\end{center}

\end{thm}

In this case we find the following subsystem of linear
inequalities.

\begin{align*}
x_1 + y_1 + z_1 & \leq x_2 + y_2 + z_2 + x_3 + y_3 + z_3 \qquad & (6,0,0)\\
x_1 + y_1 + z_2 & \leq z_1 + x_2 + y_2 + x_3 + y_3 + z_3 \qquad & (5,1,0)\\
x_1 + z_1 + y_2 & \leq y_1 + x_2 + z_2 + x_3 + y_3 + z_3 \qquad & (4,2,0)\\
x_1 + y_2 + z_2 & \leq y_1 + z_1 + x_2 + x_3 + y_3 + z_3 \qquad &
(3^{\prime},3^{\prime \prime},0)\\
\ \\
(**)\qquad x_1 + z_1 + z_2 & + z_3 \leq y_1 + x_2 + y_2 + x_3 +
y_3 \qquad &
(4,1,1)\\
x_1 + y_2 + z_3 & \leq y_1 + z_1 + x_2 + z_2 + x_3 + y_3 \qquad &
(3^{\prime},2,1)\\
(**)\qquad y_1 + z_1 + y_2 & +z_3 \leq x_1 + x_2 + z_2 + x_3 + y_3
\qquad &
(3^{\prime \prime},2,1)\\
\end{align*}

{\it After symmetrizing there are $36$ inequalities. None are
trivially redundant. However the $9$ inequalities corresponding to
the ordered partitions $(4,1,1)$ and $(3^{\prime \prime},2,1)$
and marked by ($**$) above do not occur for $\text{Sp}(6)$ and
consequently they must be redundant.}

 We now check this directly.

In order to see that the $3$ inequalities corresponding to
$(4,1,1)$ are redundant, we observe (from the first subsystem
correspoding to $G/P_1$) that $x_1 \leq x_2 + x_3$. Furthermore,
we have the inequalities for $\Delta$ given by $z_i \leq y_i,
1\leq i \leq 3$. Hence $x_1 + z_1 + z_2 +z_3 \leq x_2 + x_3 + y_1
+ y_2 + y_3$. As for the $6$ inequalities corresponding to
$(3^{\prime\prime},2,1)$, we have (from the first subsystem)
$z_1\leq z_2 + x_3$ and the inequalities (for $\Delta$) $y_1 \leq
x_1$, $y_2 \leq x_2$, and $z_3 \leq y_3$. Hence $z_1 + y_1 + y_2 +
z_3 \leq z_2 + x_3 + x_1 + x_2 + y_3$.

{\it To summarize,  for $B_3$, there are altogether  $135 = 126 +
9$ inequalities (including $9$ needed to define $\Delta^3$ in
$\mathfrak{a}^3$) of which $24$ are trivially redundant and there
are $9$ more redundant inequalities. These $9$ inequalities do no
occur in the system for $\text{\rm Sp}(6)$.
 Hence the  subsystem
for $B_3$ can be brought down to  altogether $102$ inequalities.
Moreover, a computer calculation shows that the polyhedral cone
$D_3(B_3)$ has exactly $102$ faces and thus these $102$
inequalities are irredundant.} (Of course, by Theorem
\ref{transfer}, $D_3(B_3)=D_3(C_3)$.)

\section{Generators of the cone}
In the previous section we have described the irredundant system
of linear inequalities defining the polyhedral cones $D_3(A_3)$
and $D_3(C_3) = D_3(B_3)$. We now give a system of  generators for
the cone  $D_3(C_3)$.  The components of each of the $51$
generators are arranged in the order
$x_1,x_2,x_3,y_1,y_2,y_3,z_1,z_2,z_3$ whereas the coordinates of
the corresponding generators $(v_1,v_2,v_3) \in D_3(C_3)$ satisfy
$$v_i = (x_i,y_i,z_i), 1 \leq i \leq 3.$$

\begin{thm}
The following $51$ vectors are a set of generators of the
polyhedral cone $D_3(C_3) = D_3(B_3)$ in the $9$ dimensional space
$\mathfrak{a}^3$.
\begin{center}
\setlength{\tabcolsep}{ .75cm}
\renewcommand{\arraystretch}{1.45}

\begin{tabular}{l c l c}

(1) \ 1 1 2 1 1 0 0 0 0 &  (27) \ 3 2 1 1 2 1 1 0 1 \\
(2) \ 1 1 2 1 1 0 1 1 0 &  (28) \ 1 2 2 1 1 2 0 1 0 \\
(3) \ 1 1 2 1 1 1 0 0 1 &  (29) \ 1 2 2 1 2 1 0 0 1 \\
(4) \ 1 1 2 1 1 1 0 1 0 &  (30) \ 2 1 2 1 1 2 1 0 0 \\
(5) \ 1 1 2 1 1 1 1 0 0 &  (31) \ 2 1 2 2 1 1 0 0 1 \\
(6) \ 1 2 1 1 0 1 0 0 0 &  (32) \ 2 2 1 1 2 1 1 0 0 \\
(7) \ 1 2 1 1 0 1 1 0 1 &  (33) \ 2 2 1 2 1 1 0 1 0 \\
(8) \ 1 2 1 1 1 1 0 0 1 &  (34) \ 0 1 1 0 1 1 0 0 0 \\
(9) \ 1 2 1 1 1 1 0 1 0 &  (35) \ 1 0 1 1 0 1 0 0 0 \\
(10) \ 1 2 1 1 1 1 1 0 0 & (36) \ 1 1 0 1 1 0 0 0 0 \\
(11) \ 2 1 1 0 1 1 0 0 0 & (37) \ 1 1 1 0 0 1 0 0 0 \\
(12) \ 2 1 1 0 1 1 0 1 1 & (38) \ 1 1 1 0 1 0 0 0 0 \\
(13) \ 2 1 1 1 1 1 0 0 1 & (39) \ 1 1 1 1 0 0 0 0 0 \\
(14) \ 2 1 1 1 1 1 0 1 0 & (40) \ 0 1 1 0 1 1 0 1 1 \\
(15) \ 2 1 1 1 1 1 1 0 0 & (41) \ 1 0 1 1 0 1 1 0 1 \\
(16) \ 0 1 1 0 0 0 0 0 0 & (42) \ 1 1 0 1 1 0 1 1 0 \\
(17) \ 1 1 2 1 1 2 1 1 0 & (43) \ 1 1 1 0 1 1 0 0 1 \\
(18) \ 1 2 1 1 2 1 1 0 1 & (44) \ 1 1 1 0 1 1 0 1 0\\
(19) \ 1 2 3 1 2 1 1 0 1 & (45) \ 1 1 1 1 0 1 0 0 1\\
(20) \ 1 3 2 1 1 2 1 1 0 & (46) \ 1 1 1 1 0 1 1 0 0\\
(21) \ 1 0 1 0 0 0 0 0 0 & (47) \ 1 1 1 1 1 0 0 1 0\\
(22) \ 2 1 1 2 1 1 0 1 1 & (48) \ 1 1 1 1 1 0 1 0 0\\
(23) \ 2 1 3 2 1 1 0 1 1 & (49) \ 1 1 1 1 1 1 0 0 1\\
(24) \ 1 1 0 0 0 0 0 0 0 & (50) \ 1 1 1 1 1 1 0 1 0\\
(25) \ 2 3 1 2 1 1 0 1 1 & (51) \ 1 1 1 1 1 1 1 0 0\\
(26) \ 3 1 2 1 1 2 1 1 0 &
\end{tabular}
\end{center}

\end{thm}

%\vskip5ex

%Addresses:

%S.K.: Department of Mathematics, University of North Carolina,
%Chapel Hill, NC 27599-3250, USA.\,\, (shrawan$@$email.unc.edu)

%J.M.: Department of Mathematics, University of Maryland, College
%Park, MD 20742, USA. \,\,(jjm$@$math.umd.edu)

\end{document}